\documentclass{amsart}
\usepackage{amssymb,curves,eepic}
\usepackage[dvips]{epsfig}
\usepackage[dvips]{color}
\newcommand\aut{\operatorname{Aut}}
\newcommand\C{{\mathbb C}}
\newcommand\R{{\mathbb R}}
\newcommand\Z{{\mathbb Z}}
\newcommand\N{{\mathbb N}}

\newcommand\comm{\operatorname{comm}}

\newcommand\doverline[1]{\overline{\overline{#1}}}
\newcommand\gf{{\mathcal G}}
\newcommand\hilb{{\mathcal H}}
\newcommand\msmash[1]{\raisebox{-1ex}{\makebox(0,0)[b]{\smash{$#1$}}}}
\newcommand\rist{\operatorname{rist}}
\newcommand\st{\operatorname{st}}
\newcommand\sym[1]{{\mathfrak S_{#1}}}

\newcommand\traut{{\mathcal A}}
\newcommand\tree{\mathcal T}
\newcommand\sch{{\mathcal S}}
\newcommand\spec{\operatorname{spec}}
\newcommand\uhilb{\mathcal U(\hilb)}

\newcommand\red[1]{\textcolor{red}{#1}}
\newcommand\blue[1]{\textcolor{blue}{#1}}

\DeclareMathSymbol{\idmat}{\mathbin}{AMSb}{"62}

\newtheorem{thm}{Theorem}[section]
\newtheorem{prop}[thm]{Proposition}
\newtheorem{lem}[thm]{Lemma}
\newtheorem{cor}[thm]{Corollary}
\newtheorem{defn}[thm]{Definition}
\newtheorem{scholium}[thm]{Scholium}
\newtheorem{question}{Question}
\newcommand\emdef[1]{\emph{#1}
  }

\begin{document}
\title{On the Spectrum of Hecke Type Operators related to some Fractal Groups}
\author{Laurent Bartholdi\and Rostislav I. Grigorchuk}
\date\today
\email{Laurent.Bartholdi@math.unige.ch\qquad grigorch@mi.ras.ru}
\address{\parbox{.4\linewidth}{Section de Math\'ematiques\\
    Universit\'e de Gen\`eve\\
    CP 240, 1211 Gen\`eve 24\\
    Switzerland}
  \qquad\parbox{.4\linewidth}{Steklov Institute of Mathematics\\
    Gubkina 8\\
    Moscow 117966\\
    Russia}}
\thanks{The second author wishes to express his thanks to the ``Swiss
  National Science Foundation'' and the Max-Planck Institute in Bonn
  for their support}
\begin{abstract}
  We give the first example of a connected 4-regular graph whose
  Laplace operator's spectrum is a Cantor set, as well as several
  other computations of spectra following a common ``finite
  approximation'' method. These spectra are simple transforms of the
  Julia sets associated to some quadratic maps. The graphs involved
  are Schreier graphs of fractal groups of intermediate growth, and
  are also ``substitutional graphs''. We also formulate our results in
  terms of Hecke type operators related to some irreducible quasi-regular
  representations of fractal groups and in terms of the Markovian
  operator associated to noncommutative dynamical systems via which
  these fractal groups were originally defined
  in~\cite{grigorchuk:burnside-orig}.

  In the computations we performed, the self-similarity of the groups
  is reflected in the self-similarity of some operators; they are
  approximated by finite counterparts whose spectrum is computed by an
  ad hoc factorization process.
\end{abstract}
\maketitle

\section{Introduction}
The Hecke, Markov, and Laplace operators occur in various guises
throughout mathematics. We start by a review of their more common
appearances.

\subsection{Discrete Laplacian and Hecke type Operators}
Let $\gf=(V,E)$ be a locally finite graph: there are maps
$\alpha,\omega:E\to V$ giving the extremities of edges, and every
vertex $v\in V$ has finite degree $\deg v=|\{e\in E|\,\alpha(e)=v\}|$.
Therefore all edges are oriented, and $\gf$ may have loops
($\alpha(e)=\omega(e)$) and multiple edges. The \emdef{discrete
  Laplace operator} of $\gf$ is the operator $\Delta=1-M$ on
$\ell^2(V,\deg)$, where $M$ is the ``adjacency'' or Markovian operator
\[(Mf)(v)=\frac1{\deg v}\sum_{e\in E:\,\alpha(e)=v}f(\omega(e)).\]
The
theory of discrete Laplace operators $\Delta$ has a long history, and
is a popular topic of contemporary
mathematics~\cite{cvetkovic:spectra,woess:cogrowth}. In the context of
random walks on graphs, one usually considers the Markovian operator
$M$ rather than $\Delta$: if $e_v$ be the Dirac delta function at the
vertex $v$, then $\langle M^n e_v|e_w\rangle$ is the probability of a
random walk starting at $v$ to reach $w$ in $n$ steps. If the random
walk is symmetric, then $M$ is a self-adjoint operator and its
spectrum lies in $[-1,1]$. The spectral properties of $M$ contain
valuable information for the theory of random walks and discrete
potential, graph theory, abstract harmonic analysis, the theory of
operator algebras, etc. For instance, a theorem by Harry
Kesten~\cite{kesten:rwalks}, generalized by different mathematicians
(see~\cite{woess:cogrowth} and~\cite{ceccherini-:amen} with its
bibliography) asserts that $\gf$ is amenable if and only if $1$ is in
the spectrum of $M$. Note that the random walk need not be simple
(probabilities of moving in different directions may be different); a
Markovian operator can still be associated to the walk.  Let us
finally mention a more general setting, developed these last
years~\cite{novikov:schrodinger}: the (discrete) Schr\"odinger
operators $\Delta+P$, where $P$ is diagonal, and the coefficients of
$\Delta$ may depend on the vertex they correspond to.

The theory of Hecke type operators was developed in parallel: if
$\pi:G\to\uhilb$ is a unitary representation of a finitely generated
group $G$ given with a symmetric generating system
$S=\{s_1,\dots,s_m\}=S^{-1}$ in a Hilbert space $\hilb$, then one
associates to $\pi$ a self-adjoint operator $H$ on $\hilb$:
\[H = \sum_{i=1}^m p(i)\pi(s_i),\]
for some $p(i)\in\C$. The most important choice is $p(i)=\frac1m$ for
all $i\in\{1,\dots,m\}$; we shall restrict to this choice in the
sequel, and assume, when no weight is given, that this one is used.

Operators of Hecke type play an important role in mathematical
physics~\cite{connes:ncg}, Arakelov theory in number theory
(see~\cite{li:nt},~\cite{serre:hecke} and~\cite{serre:arith} for the
connection between number theory and operators), and Ramanujan
graphs~\cite{lubotzky:expanders}.  The group-theoretical content of
$H$, mainly in the case of the regular representation, was studied by
Pierre de la Harpe, A. Guyan Robertson and Alain Valette
in~\cite{harpe-r-v:sg,harpe-r-v:sg2}, and in many other papers --- see
the bibliography in~\cite{woess:cogrowth}.

\subsection{Spectra of Noncommutative Dynamical Systems}
Let $T$ be an invertible measure-preserving transformation of a
measure space $(X,\mu)$, and let $A$ be the corresponding unitary
operator in $L^2(X,\mu)$, given by $(Af)(x)=f(T^{-1}x)$. By the
\emdef{spectrum} of the dynamical system one usually means the
spectrum of the operator $A$, or, as is almost the same, the spectrum
of the self-adjoint operator $A+A^{-1}$.  These last spectra are in
correspondence through the map $z\mapsto z+z^{-1}$. If $T$ is
aperiodic, then $\spec(A+A^{-1})=[-1,1]$.

By a \emdef{noncommutative dynamical system} we mean a collection $S$
of invertible measure-class-preserving transformations on a measure space
$(X,\mu)$, that do not necessarily commute. Let $G$ be the group
generated by $S$. It has a natural unitary representation $\pi$ in
$\mathcal(L^2(X,\mu))$ given by
\[(\pi(g)f)(x) = \sqrt{\mathfrak g(x)}f(g^{-1}x),\]
where $\mathfrak g(x)=dg\mu(x)/d\mu(x)$ is the Radon-Nikod\'ym
derivative.  The \emdef{spectrum} of the dynamical system $S$ is the
spectrum of the Hecke type operator associated to $G$, $S\cup S^{-1}$ and
$\pi$. In case $|S|=1$, this definition reduces to the previous,
classical one.

\subsection{Examples} One of the most famous operators of
Hecke type is the Harper-Mathieu-Peierls operator $H_\lambda$ on
$\ell^2(\Z)$ acting on infinite sequences $f:\Z\to\C$ by
\[\left(H_\lambda(f)\right)(n)= f(n-1)+f(n+1)+(2\cos\lambda n)f(n),\]
for any $\lambda\in\R$. These $H_\lambda$ are the Hecke type operators
associated to the Heisenberg group
\[\left.\left\{\begin{pmatrix}1&m&p\\&1&n\\&&1\end{pmatrix}\right|\,m,n,p\in\Z\right\}\cong\langle a,b,c|\,[a,b]=c,[a,c]=[b,c]=1\rangle\]
and its representation $\pi_\lambda$ in $\ell^2(\Z)$, where
$\pi_\lambda(a)$ acts by translation: $[\pi_\lambda(a)](f)(n)=f(n-1)$,
and $\pi_\lambda(b)$ acts by pointwise multiplication with the
function $e^{i\lambda n}$, namely $[\pi_\lambda(b)](f)(n)=e^{i\lambda
  n}f(n)$.

The Harper operator is the operator related to the Quantum Hall
effect, and originally arose in connection with the two-dimensional
lattice. One can start from any Cayley graph, and construct a
corresponding Harper operator, which would be the discrete analogue of
the magnetic Laplacian~\cite{carey-:hall}.

The spectral properties of this operator were thoroughly investigated;
if $\lambda$ is a Liouville number, the spectrum of $H_\lambda$ is a
Cantor set~\cite{bellissard-s:cantor}. We note that $H_\lambda$ is a
Schr\"odinger operator on the one-dimensional lattice. By Fourier
transform, it can be realized as an element of the cross product
$C^*$-algebra $R_\lambda\ltimes\mathcal C(\mathbb S^1)$, where
$R_\lambda$ is the dynamical system generated by an angle-$\lambda$
rotation on $\mathbb S^1$ and $\mathcal C(\mathbb S^1)$ denotes the
algebra of continuous functions on the circle.

Another example was studied by David Kazhdan~\cite{kazhdan:uniform}.
Let $\alpha$ and $\beta$ be two noncommuting rotations in the plane
$\R^2$. They generate a group $G$ with a unitary action $\pi$ on
$L^2(\R^2)$. Kazhdan studies the operator
$M=\pi(\alpha)+\pi(\alpha^{-1})+\pi(\beta)+\pi(\beta^{-1})$ and shows
that its Fourier transform decomposes as a direct integral of
operators acting on functions on the circle. Each of these is an
element of $R_\lambda\ltimes\mathcal B(\mathbb S^1)$, where $\mathcal
B(\mathbb S^1)$ denotes the bounded functions on the circle, and
happens to be a Schr\"odinger operator; spectral properties of these
operators are then used to show that the orbits of $G$ in $\R^2$ are
uniformly distributed.

\subsection{Main Results} We produce examples of operators of Hecke
type with Cantor set spectrum, but where additionally the
representations $\pi$ involved are quasi-regular. This produces graphs
whose Laplace operators have totally discontinuous spectrum, namely the
associated Schreier graphs. Our main results read:
\begin{thm}\label{thm:main}
  \begin{enumerate}
  \item There is a connected $4$-regular graph of polynomial growth,
    which is a Schreier graph of a group of intermediate growth, and
    whose Laplacian's spectrum is a Cantor set.
  \item There is a connected $4$-regular graph of polynomial growth,
    which is a Schreier graph of a group of intermediate growth, and
    whose Laplacian's spectrum is the union of a Cantor set $K$ and a
    countable set $P$ of isolated points whose accumulation set is
    $K$.
  \item There are noncommutative dynamical systems generated by $2$
    transformations whose spectrum are the same as in the above two
    points.
  \item The above spectra are calculated explicitly. The Cantor set
    $K$ is of the form $F(J)$ where $F$ is a simple algebraic function
    and $J$ is the Julia set of a quadratic map $z\mapsto
    z^2-\lambda$, where $\lambda=45/16$ in the first case and
    $\lambda=6$ in the second case. $J$ is the set of points of the
    form
    \[\pm\sqrt{\lambda\pm\sqrt{\lambda\pm\sqrt{\lambda\pm\sqrt{\dots}}}}.\]
  \end{enumerate}
\end{thm}
To the best of our knowledge, these are the first examples of graphs
of constant vertex degree whose spectrum is totally disconnected.
There are, however, examples of Schr\"odinger operators on $\Z$ (or
$\R$) with nowhere dense spectrum; they are obtained as Harper
operators (as mentioned above) or following a result by J\"urgen
Moser~\cite{moser:schrodinger}.

There are also examples of random walks on non-regular graphs, but
with vertex degrees $1$ or $3$, whose spectrum is the union of a
countable set and a Cantor set of null Lebesgue
measure~\cite{malozemov:fractal}.

Similarly we produce Hecke type operators of quasi-regular representations
that have the same spectra as above. These are probably the first
examples of quasi-regular representations of virtually torsion-free
groups with totally discontinuous spectrum; at least, the
representations $\pi_\lambda$ of the Heisenberg group are not
quasi-regular, but come from the cross-product construction, which is
often used to produce interesting examples. In the sequel we produce
interesting examples of spectra using purely non-commutative dynamical
systems and associated methods.

The graphs mentioned in Theorem~\ref{thm:main} are Schreier graphs of
some fractal groups. They are of polynomial growth and have a clear
``fractal'' appearance; see Figure~\ref{fig:schreier}. By a fractal group
we mean a group which acts on a regular rooted tree $\tree$, such that
this action has some self-similar properties; this notion is very much
related to that of branch group introduced
in~\cite{grigorchuk:jibg}. The Schreier graphs $\sch(G,P,S)$ are
defined in~\ref{defn:sch}; in our examples we take for $P$ the
stabilizer $\st_G(e)$ of an infinite ray starting at the root of
$\tree$, i.e.\ an element of the boundary $\partial\tree$.

The parabolic subgroups $P=\st_G(e)$ have the remarkable property of
being \emph{weakly maximal}: $[G:P]=\infty$ but $[G:L]<\infty$ for all
$L\gneqq P$. The corresponding quasi-regular representations
$\rho_{G/P}$ are irreducible, and $\bigcap_{g\in G}P^g=1$. We thus
have an important family of faithful irreducible representations of
$G$, that deserves further investigation.

The first example of group of fractal type was constructed
in~\cite{grigorchuk:burnside-orig} as an example of infinite torsion
$2$-group; it was described as a set of measure-preserving
transformations of the interval $[0,1]$, but can equivalently be
described by its action on a rooted tree (the binary expansion of a
real in $[0,1]$ giving a path in the rooted binary tree). Later many
new examples of this sort
appeared~\cite{gupta-s:infinitep,grigorchuk:growth,bartholdi-g:parabolic}.
It then became clear that the study of these groups via their tree
action was most fruitful and led to interesting ideas and results; for
a survey see~\cite{grigorchuk:jibg}; and for an introduction to the
first example, $G$, see~\cite{harpe:cgt,ceccherini-:grigp}.

Among the five groups we consider, two ($\Gamma$ and
$\overline\Gamma$) are virtually torsion-free, and have a totally
disconnected spectrum. The existence of such groups lends some hope to
the existence of torsion-free fractal groups whose Laplace operator
has a totally disconnected spectrum, or at least a gap in the
spectrum. Such an example would provide a counterexample to the
Kaplansky-Kadison conjecture on idempotents (that implies that the
spectrum of the Laplace operator related to the regular representation
is connected), and to the Baum-Connes
conjecture~\cite{valette:idempotents}. Note that if such a group
existed, it would be non-amenable~\cite{higson-k:bc}.

The method used in the computations is the following: we compute
explicitly the spectrum of the finite graph $\sch(G,P_n,S)$, where
$P_n$ is the stabilizer of the rightmost vertex in the $n$th row of
the tree $\tree$. Then some arguments, coming
from~\cite{lubotzky:cayley,grigorchuk-z:infinite,martin-v:bs} but
adapted to our goals, are used to obtain the spectrum of the infinite
graphs from the finite spectra.

\subsection{Notation} We assume all groups act on the left on sets,
and write $g^h=hgh^{-1}$ and $[g,h]=ghg^{-1}h^{-1}$. We also write
$\langle S\rangle$ and $\langle S\rangle^G$ for the subgroup and
normal subgroup of $G$ generated by $S$. The regular representation of
$G$ in $\ell^2(G)$ is written $\rho_G$, and the quasi-regular
representation of $G$ in $\ell^2(G/H)$ is written $\rho_{G/H}$.
The symmetric group on a set $S$ of cardinality $n$ is written $\sym
S$ or $\sym n$.

\subsection{Guide to Quick Reading}
We present in this paper five computations of spectra related to
groups; however, the reader interested solely in examples of graphs
with Cantor-set spectrum may wish to skip the group-theoretic
discussion. In this case the sections~\ref{subs:computeGamma}
and~\ref{sec:schreier} should describe the construction in a fairly
self-contained manner.

\section{Groups acting on rooted trees}
The groups we shall consider will all be subgroups of the group
$\aut(\tree)$ of automorphisms of a regular rooted tree $\tree$. Let
$\Sigma$ be a finite alphabet.  The vertex set of the tree
$\tree_\Sigma$ is the set of finite sequences over $\Sigma$; two
sequences are connected by an edge when one can be obtained from the
other by right-adjunction of a letter in $\Sigma$.  The top node is
the empty sequence $\emptyset$, and the children of $\sigma$ are all
the $\sigma s$, for $s\in\Sigma$. We suppose $\Sigma=\Z/d\Z$, with the
operation $\overline s=s+1\mod d$.  Let $a$, called the \emdef{rooted
  automorphism} of $\tree_\Sigma$, be the automorphism of
$\tree_\Sigma$ defined by $a(s\sigma)=\overline s\sigma$: it acts
nontrivially on the first symbol only, and geometrically is realized
as a cyclic permutation of the $d$ subtrees just below the root.
\begin{center}
\epsfig{file=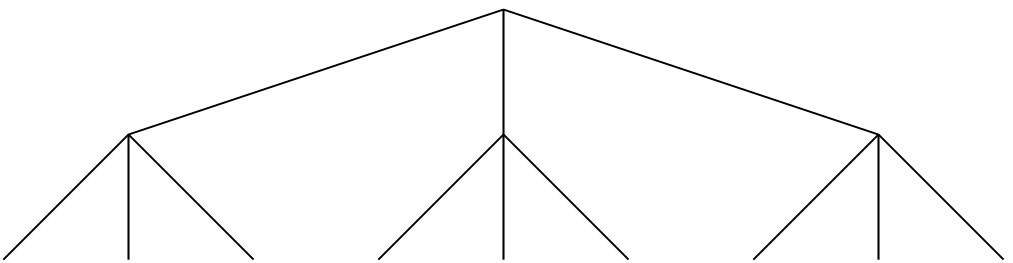}
\end{center}
Fix some $\Sigma$ and let $\traut=\aut(\tree_\Sigma)$. For any
subgroup $G<\traut$, let $\st_G(\sigma)$ denote the subgroup of $G$
consisting of the automorphisms that fix the sequence $\sigma$, and
$\st_G(n)$ denote the subgroup of $G$ consisting of the automorphisms
that fix all sequences of length $n$:
\[\st_G(\sigma)=\{g\in
G|\,g\sigma=\sigma\},\qquad\st_G(n)=\bigcap_{\sigma\in\Sigma^n}\st_G(\sigma).\]
The $\st_G(n)$ are normal subgroups of finite index of $G$; in
particular $\st_G(1)$ is of index at most $d!$. Let $G_n$ be the
quotient $G/\st_G(n)$. If $g\in\traut$ is an automorphism fixing the
sequence $\sigma$, we denote by $g_{|\sigma}$ the element of $\traut$
corresponding to the restriction to sequences starting by $\sigma$:
\[\sigma g_{|\sigma}(\tau)=g(\sigma\tau).\]
As the subtree starting from any vertex is isomorphic to the initial
tree $\tree_\Sigma$, we obtain this way a map
\begin{equation}\label{eq:phi}
  \phi:\begin{cases}\st_\traut(1)\to\traut^\Sigma\\h\mapsto (h_{|0},\dots,h_{|d-1})\end{cases}
\end{equation}
which is an embedding.

\begin{defn}
  A subgroup $G<\traut$ is \emph{level-transitive} if the action of
  $G$ on $\Sigma^n$ is transitive for all $n\in\N$. We shall always
  implicitly make that assumption.

  $G$ is \emdef{fractal} if for every vertex $\sigma$ of
  $\tree_\Sigma$ one has $\st_G(\sigma)_{|\sigma}\cong G$, where the
  isomorphism is given by identification of $\tree_\Sigma$ with its
  subtree rooted at $\sigma$.
\end{defn}

For a sequence $\sigma$ and an automorphism $g\in\traut$, we
denote by $g^\sigma$ the element of $\traut$ acting as $g$ on
the sequences starting by $\sigma$, and trivially on the others:
\[g^\sigma(\sigma\tau)=\sigma g(\tau),\qquad g^\sigma(\tau)=\tau\text{
  if }\tau\text{ doesn't start by }\sigma.\]
Let $G<\traut$ be a group acting faithfully, and transitively on each level,
on a rooted tree $\tree_\Sigma$. The \emdef{rigid stabilizer} of $\sigma$ is
$\rist_G(\sigma)=\{g^\sigma|\,g\in G\}\cap G$. We say $G$ has
\emdef{infinite rigid stabilizers} if all the $\rist_G(\sigma)$ are
infinite.

\begin{defn}\label{defn:b}\begin{enumerate}
  \item $G$ is a \emdef{regular branch} group if it has a finite-index
    subgroup $K<\st_G(1)$ such that
    \[K^\Sigma<\phi(K).\]\label{defn::rb}
  \item A subgroup $G<\traut$ is a \emdef{branch group} if for every
    $n\ge1$ there exists a subgroup $L_n<\traut$ and an embedding
    \[L_n\times\dots\times L_n\hookrightarrow \st_G(n),\]
    where the direct product is indexed by $\Sigma^n$, the injection
    is given on each factor by $(\ell,\sigma)\mapsto\ell^\sigma$, and
    the image is normal of finite index in $\st_G(n)$.\label{defn::b}
  \item $G$ is a \emdef{weak branch} group if all of its rigid
    stabilizers $\rist_G(\sigma)$ are infinite.\label{defn::wb}
  \end{enumerate}
\end{defn}
If $G$ is fractal, one has for all $n$ an embedding
$\st_G(n)<G^{d^n}$. Note that the definition of a branch group admits
an even more general setting --- see~\cite{grigorchuk:jibg}. Four of
our examples will be regular branch groups, and the last one will not
be a branch, but rather a weak branch group. The following lemma shows
that, for fractal groups, \ref{defn::rb} implies \ref{defn::b} implies
\ref{defn::wb} in Definition~\ref{defn:b}.
\begin{lem}\label{lem:rb->b}
  If $G$ is a fractal, regular branch group, then it is a branch group.
  If $G$ is a branch group, then it is a weak branch group.
\end{lem}
\begin{proof}
  Assume $G$ is a regular branch group on its subgroup $K$. Define
  $L_n=K$ for all $n$. Clearly $\st_G(n)$ contains the direct product
  $L_n^{\Sigma^n}$, and it is of finite index in $G^{\Sigma^n}$, so
  all the more of finite index in $\st_G(n)$.  The second implication
  holds because branch groups are infinite, and `finite index in
  infinite group' is stronger than `infinite'.
\end{proof}

In the sequel we shall be concerned with subgroups $G$ of $\traut$
that are finitely generated, fractal, and contain the rooted
automorphism $a$. These groups will be naturally equipped with the
restriction of the map $\phi$ defined in~(\ref{eq:phi}), a descending
sequence of normal subgroups $H_n=\st_G(n)$ and an approximating
sequence of finite quotients $G_n=G/H_n$. These quotients can be seen
as subgroups of the symmetric group $\sym{\Sigma^n}$ on
$\Sigma^n$. More details on all of these groups and their subgroups
appear in~\cite{bartholdi-g:parabolic}.

\subsection{Dynamical Systems}\label{subs:ds}
Assume as above that a group $G$ generated by a set $S$ acts on the
$d$-regular rooted tree $\tree=\{0,\dots,d-1\}^*$. Then $G$ acts
naturally on the boundary $\partial\tree=\{0,\dots,d-1\}^\N$, and this
action preserves the uniform Bernoulli measure $\nu$ on the compact
space $\partial\tree$. We associate thus a dynamical system
$(G,S,\partial\tree,\nu)$ to the group $G$.

This dynamical system is naturally isomorphic to a dynamical system
$(G,S,[0,1],m)$, where $m$ is the Lebesgue measure, and $G$ (generated
by $S$) acts on $[0,1]$ by measure-preserving transformations in the
following way: let $g\in G$, and $\gamma\in[0,1]$ a $d$-adic
irrational with base-$d$ expansion $0.\gamma_1\gamma_2\dots$. Then
$g(\gamma)=0.\delta_1\delta_2\dots$, where the infinite sequence
$(\gamma_1,\gamma_2,\dots)$ is mapped by $g$ to
$(\delta_1,\delta_2,\dots)$.  This defines the action of $G$ on a
subset of full measure of $[0,1]$.

The orbits of $G$ on $\partial\tree$ can be made explicit as follows:
\begin{defn}
  Two infinite sequences $\sigma,\tau:\N\to\Sigma$ are
  \emdef{confinal} if there is an $N\in\N$ such that $\sigma_n=\tau_n$
  for all $n\ge N$.

  Confinality is an equivalence relation, and equivalence classes are
  called \emdef{confinality classes}.
\end{defn}

\begin{prop}\label{prop:confinal}
  Let $G$ be a group acting on a regular rooted tree $\tree$, and
  assume that for any generator $g\in G$ and infinite sequence
  $\sigma$, the sequences $\sigma$ and $g\sigma$ differ only in
  finitely many places. Then the confinality classes of the action of
  $G$ on $\partial\tree$ are unions of orbits. If moreover
  $\st_G(\sigma)$ contains the rooted automorphism $a$ for all
  $\sigma\in\tree$, the orbits of the action are confinality classes.
\end{prop}
The dynamics of the actions of a group on the boundary of a tree from
the point of view of confinality are investigated in more detail in
in~\cite{nekrashevych-s:orbits}. We remark that the five example
groups---$G,\tilde G,\Gamma,\overline\Gamma,\doverline\Gamma$---we
shall consider satisfy the conditions of the proposition above.

\subsection{Growth of Groups and Parabolic Subgroups}
We recall some facts about word-growth of groups and sets on which
they act.
\begin{defn}
  Let $G$ be a group generated by a finite set $S$, let $X$ be a set
  upon which $G$ acts transitively, and choose $x\in X$. The
  \emdef{growth} of $X$ is the function $\gamma:\N\to\N$ defined by
  \[\gamma(n)=|\{gx\in X|\,|g|\le n\}|,\]
  where $|g|$ denotes the minimal length of $g$ when written as a word
  over $S$.  By the \emdef{growth of $G$} we mean the growth of the
  action of $G$ on itself by left-multiplication.
  
  Given two functions $f,g:\N\to\N$, we write $f\preceq g$ if there is
  a constant $C\in\N$ such that $f(n)<Cg(Cn+C)+C$ for all $n\in\N$,
  and $f\sim g$ if $f\preceq g$ and $g\preceq f$. The equivalence
  class of the growth of $X$ is independent of the choice of $S$ and
  of $x$.

  $X$ is of \emdef{polynomial growth} if $\gamma(n)\preceq n^d$ for some
  $d$. It is of \emdef{exponential growth} if $\gamma(n)\succeq
  e^n$. It is of \emdef{intermediate growth} in the remaining
  cases. This trichotomy does not depend on the choice of $x$.
\end{defn}

Assume now that $G$ is a group acting on the tree $\Sigma^*$,
and that a subset $S\subset G$ is given.
\begin{defn}
  The \emdef{portrait} of $g\in G$ with respect to $S$ is a subtree of
  $\Sigma^*$, with inner vertices labeled by $\sym\Sigma$ and leaf
  vertices labeled by $S\cup\{1\}$. It is defined recursively as follows: if
  $g\in S\cup\{1\}$, the portrait of $g$ is the subtree reduced to the root
  vertex, labeled by $g$ itself. Otherwise, let $\alpha\in\sym\Sigma$
  be the permutation of the top branches of $\Sigma^*$ such that
  $g\alpha^{-1}\in\st_G(1)$; let
  $(g_0,\dots,g_{d-1})=\phi(g\alpha^{-1})$ and let $\tree_i$ be the
  portrait of $g_i$. Then the portrait of $g$ is the subtree of
  $\Sigma^*$ with $\alpha$ labeling the root vertex and subtrees
  $\tree_0,\dots,\tree_{d-1}$ connected to the root.

  The \emdef{depth} of $g\in G$ is the height (length of a maximal
  path starting at the root vertex) $\partial(g)\in\N\cup\{\infty\}$
  of the portrait of $g$.
\end{defn}
Therefore the depth of $g$ is finite if and only if the portrait of
$G$ is finite. Both are finite for all groups we consider in this
paper, and the depth may be estimated using the following lemma:
\begin{lem}\label{lem:finitep}
  Assume $S$ generates $G$ and $\phi:g\mapsto(g_1,\dots,g_d)$ defined
  in~(\ref{eq:phi}) has the property that $|g_i|<|g|$ for all $i$.
  Then every $g\in G$ has a finite portrait. If moreover there are
  constants $\lambda,\mu$ with $\mu/(1-\lambda)<2$ such that
  $|g_i|\le\lambda|g|+\mu$ for all $i$, then asymptotically when
  $|g|\to\infty$
  \[\partial(g)\le\log_{1/\lambda}|g|.\]
\end{lem}
\begin{proof}
  Consider the `level' function
  \[F(n)=\max_{|g|<n}\partial(g).\]
  Then $F$ is increasing, and by assumption $F(2)=0$, $F(n)\le F(\lambda
  n+\mu)+1$. It then follows that
  \[F(n) \le F(\lambda n+\mu)+1\le\dots\le F(\lambda^k n+\lambda^{k-1}\mu+\dots+\lambda\mu+\mu)+k;\]
  let us take for $k$  a natural number satisfying
  $\lambda^k n+\lambda^{k-1}\mu+\dots+\lambda\mu+\mu\le 2$, for
  instance
  \[k=\left\lceil\log_\lambda\frac{2-\mu/(1-\lambda)}{n-mu/(1-\lambda)}\right\rceil,\]
  where $\lceil x\rceil$ is the least integer greater than $x$. The
  result follows, because then $F(n)\le k$ and $k\sim \log_{1/\lambda}n$.
\end{proof}

\begin{scholium}
  For all groups considered in this paper, we have
  $|g_i|\le\frac12|g|+\frac12$ for their natural generating systems,
  as can be checked on the tables describing $\phi$. As a consequence,
  they all satisfy $\partial(g)\precsim\log_2|g|$.
\end{scholium}
 
\begin{defn}
  Let $\tree=\Sigma^*$ be a rooted tree. A \emdef{ray} $e$ in $\tree$
  is an infinite geodesic starting at the root of $\tree$, or
  equivalently an element of $\partial\tree=\Sigma^\N$.
  
  Let $G<\traut$ and $e$ be a ray. The associated \emdef{parabolic
    subgroup} is $\st_G(e)=\cap_{n\ge0}\st_G(e_n)$, where $e_n$ is the
  length-$n$ prefix of $e$.
\end{defn}
Assume that $G$ satisfies the conditions of Lemma~\ref{lem:finitep}.
Then we have the
\begin{prop}\label{prop:polygrowth}
  Let $G<\traut$ satisfy the conditions of
  Proposition~\ref{prop:confinal} and Lemma~\ref{lem:finitep} (for the
  constant $\lambda$), and let $P$ be a parabolic subgroup. Then
  $G/P$, as a $G$-set, is of polynomial growth of degree at most
  $\log_{1/\lambda}(d)$.
  If moreover $G$ is level-transitive, then $G/P$'s asymptotical
  growth is polynomial of degree $\log_{1/\lambda}(d)$.
\end{prop}
\begin{proof}
  Suppose that $P=\st_G(e)$. Then $G/P$ can be identified with the
  $G$-orbit of $e$, and, by Proposition~\ref{prop:confinal}, with the
  set of all infinite sequences over $\Sigma$ that eventually coincide
  with $e$. If $\partial(g)< k$, it sends the infinite sequence $e$ to
  one of the $d^k$ sequences in $\Sigma^ke_k\dots$; thus the image of
  $e$ under the set of elements of depth at most $k$ is of cardinality
  bounded by $d^k$.  The image of $e$ under the set of elements of
  length at most $n$ is then asymptotically bounded by
  $d^{\log_{1/\lambda}(n)}=n^{\log_{1/\lambda}(d)}$, by
  Lemma~\ref{lem:finitep}.
\end{proof}

We now recall some facts on commensurators:
\begin{defn}
  The \emdef{commensurator} of a subgroup $H$ of $G$ is
  \[\comm_G(H) = \{g\in G|\,H\cap H^g\text{ is of finite index in
    }H\text{ and }H^g\}.\]
  Equivalently, letting $H$ act on the left on the  cosets $\{gH\}$,
  \[\comm_G(H) = \{g\in G|\,H\cdot(gH)\text{ and }H\cdot(g^{-1}H)\text{ are finite orbits}\}.\]
\end{defn}

\begin{prop}[\cite{bartholdi-g:parabolic}]\label{prop:comm}
  Let $G$ be a weak branch group and let $P$ be a parabolic subgroup.
  Then $\comm_G(P)=P$.
\end{prop}
\begin{thm}[Mackey~\cite{mackey:representations,burger-h:irreducible}]\label{thm:mackey}
  Let $P<G$ be any subgroup inclusion. Then the quasi-regular
  representation $\rho_{G/P}$ is irreducible if and only if
  $\comm_G(P)=P$.
\end{thm}
Therefore, for all weak branch groups $\rho_{G/P}$ is irreducible.

The following lemma is well known:
\begin{lem}\label{lem:factorqr}
  Let $H<G$ be a subgroup of finite index. Then the orbits of $H$ on
  $G/H$, the double cosets $HgH<G$ and the irreducible
  components of the $G$-space $\ell^2(G/H)$ are all in bijection.
\end{lem}
Therefore, the orbits of $P_n=P\cdot\st_G(n)$ on $\Sigma^n$ are in
bijection with the decomposition of $\rho_{G/P_n}$ in irreducible
subrepresentations.

\subsection{Groups and Finite Automata}
There are various uses of automata in group theory, most notably as
\emph{word acceptors}, where the automata recognize some words as
group elements and perform operations on these
words~\cite{epstein-:wp}; and as \emph{transducers} or
\emph{sequential machines} (see~\cite[Chapter~XI]{eilenberg:a}
or~\cite{gecseg-c:ata}), where the automata themselves are the
elements of the group, and are distinguished by the transformation
they perform on their input. The former use gives rise to the theory
of \emph{automatic groups}; we propose to call the latter
\emph{automata groups}. The automata they are built with are called
\emph{Mealy machines} or~\emph{Moore
  machines}~(see~\cite{glushkov:ata} or~\cite[page 109]{brauer:automaten}).

We present a restricted definition of \emph{finite transducers}. In
the standard terminology, they would be called \emph{invertible
  transducers}.
\begin{defn}
  Let $\Sigma$ be a finite alphabet.  A \emph{finite transducer} on
  $\Sigma$ is a finite directed graph $\gf=(V,E)$, a labeling
  $\lambda:E\to\Sigma$ of the edges such that for each vertex $v\in V$
  the restriction of $\lambda$ is a bijection between $\{e\in
  E|\,\alpha(e)=v\}$ and $\Sigma$, and a labeling
  $\tau:V\to\sym\Sigma$ of the nodes (called \emph{states}) by the
  symmetric group on $\Sigma$.
  
  An \emph{initial transducer} $\gf_q$ is a finite transducer $\gf$
  with a distinguished initial state $q\in V$.
\end{defn}

Let $\gf$ be a finite transducer, and $\{\gf_q\}_{q\in V}$ be the set
of its initial transducers. Each $\gf_q$ defines an automorphism
$\overline{\gf_q}$ of the rooted tree $\tree_\Sigma$ as follows: let
$\sigma=\sigma_0\dots\sigma_n$ be a vertex of $\tree$. Let $e$ be the
edge of $\gf_q$ labeled $\sigma_0$. Define recursively
\[\overline{\gf_q}(\sigma_0\dots\sigma_n)=\tau(q)(\sigma_0)\overline{\gf_{\omega(e)}}(\sigma_1\dots\sigma_n).\]
We shall call two initial transducers $\gf_q$ and $\gf'_{q'}$
\emdef{equivalent} if their actions $\overline{\gf_q}$ and
$\overline{\gf'_{q'}}$ on $\Sigma^*$ are the same. Every initial
transducer is equivalent to a unique transducer that is minimal with
respect to its number of nodes~\cite[Chapter~XII, Theorem~4.1]{eilenberg:a}.

Define now $G(\gf)$ as the group generated by the $\overline{\gf_q}$,
where $q$ ranges over the set of states of $\gf$.
We call such a group an \emdef{automata group}. The following fact is
well known, and dates back to Ji\v r\'\i\ Ho\v rej\v s in the early
60's~\cite{horejs:automata}:
\begin{prop}
  Let $\gf_q$ and $\gf'_{q'}$ be finite initial transducers on the same
  alphabet $\Sigma$. Then $\overline{\gf_q}^{-1}$ and
  $\overline{\gf_q}\circ\overline{\gf'_{q'}}$ can be represented as
  finite initial transducers.
\end{prop}

In general different transducers can generate isomorphic groups. For
instance, consider the three-vertex transducer in the middle of
Figure~\ref{fig:autoGamma}. The group it generates is isomorphic (and
even conjugate in $\traut$) to the one generated by the following
two-state transducer $\gf$ on $\Sigma=\{0,1,2\}$, because the elements
$\tilde t,\tilde a$ satisfy the same recursions as $t,a$ (see
Subsection~\ref{subs:rec}):
\begin{center}\begin{picture}(150,50)(-25,-25)
    \put(0,0){\circle{20}\msmash1}\put(10,-10){\msmash{\tilde t}}
    \put(100,0){\circle{20}\msmash\varepsilon}\put(90,-10){\msmash{\tilde a}}
    \put(10,0){\vector(1,0){80}}\put(50,-5){\msmash{0,1}}
    \put(-14,0){\arc{20}{.7853981635}{5.497787145}}\put(-7,-7){\vector(1,1){0}}
    \put(-18,0){\msmash2}
    \put(114,0){\arc{20}{3.926990818}{8.639379799}}\put(107,7){\vector(-1,-1){0}}
    \put(118,0){\msmash\Sigma}
\end{picture}\end{center}
Here $\varepsilon$ is the cycle $(0,1,2)\in\sym3$. The actions
of $\gf_{\tilde t}$ and $\gf_{\tilde a}$ are as follows:
\begin{align*}
  \gf_{\tilde t}(2\dots2\sigma_m\dots\sigma_n) &=2\dots2\varepsilon(\sigma_m)\dots\varepsilon(\sigma_n)\qquad\text{when }\sigma_m\neq2,\\
  \gf_{\tilde a}(\sigma_0\dots\sigma_n) &=\varepsilon(\sigma_0)\dots\varepsilon(\sigma_n).
\end{align*}

\subsection{The Group $G$}\label{subs:defG}
We give here some basic facts about the first of our examples, the
group $G$~\cite{grigorchuk:burnside-orig,grigorchuk:gdegree-orig}. We take
$\Sigma=\{0,1\}$.  Recall $a$ is the automorphism permuting the top
two branches of $\tree_2$. Let recursively $b$ be the automorphism
acting as $a$ on the right branch and $c$ on the left, $c$ be the
automorphism acting as $a$ on the right branch and $d$ on the left,
and $d$ be the automorphism acting as $1$ on the right branch and $b$
on the left. In formul\ae,
\begin{alignat*}{2}
  b(0x\sigma)&=0\overline x\sigma,&\qquad b(1\sigma)&=1c(\sigma),\\
  c(0x\sigma)&=0\overline x\sigma,&\qquad c(1\sigma)&=1d(\sigma),\\
  d(0x\sigma)&=0x\sigma,&\qquad d(1\sigma)&=1b(\sigma).
\end{alignat*}
$G$ is the group generated by $\{a,b,c,d\}$.  It is readily checked
that these generators are of order $2$ and that $\{1,b,c,d\}$
constitutes a Klein group; one of the generators $\{b,c,d\}$ can thus be
omitted.

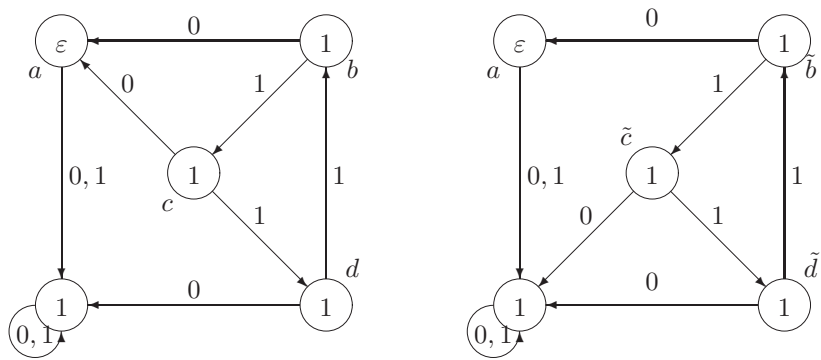
\begin{figure}
\begin{picture}(150,150)(-25,-25)
  \put(0,0){\circle{20}\msmash1}
  \put(0,100){\circle{20}\msmash\varepsilon}\put(-10,90){\msmash a}
  \put(100,0){\circle{20}\msmash1}\put(110,15){\msmash d}
  \put(100,100){\circle{20}\msmash1}\put(110,90){\msmash b}
  \put(50,50){\circle{20}\msmash1}\put(40,40){\msmash c}
  \put(90,100){\vector(-1,0){80}}\put(50,107){\msmash{0}}
  \put(0,90){\vector(0,-1){80}}\put(10,50){\msmash{0,1}}
  \put(90,0){\vector(-1,0){80}}\put(50,7){\msmash{0}}
  \put(100,10){\vector(0,1){80}}\put(105,50){\msmash{1}}
  \put(43,57){\vector(-1,1){36}}\put(25,85){\msmash{0}}
  \put(57,43){\vector(1,-1){36}}\put(75,35){\msmash{1}}
  \put(93,93){\vector(-1,-1){36}}\put(75,85){\msmash{1}}
  \put(-10,-10){\arc{20}{0}{4.712388981}}\put(0,-10){\vector(0,1){0}}\put(-10,-10){\msmash{0,1}}
\end{picture}
\qquad
\begin{picture}(150,150)(-25,-25)
  \put(0,0){\circle{20}\msmash1}
  \put(0,100){\circle{20}\msmash\varepsilon}\put(-10,90){\msmash{a}}
  \put(100,0){\circle{20}\msmash1}\put(110,15){\msmash{\tilde d}}
  \put(100,100){\circle{20}\msmash1}\put(110,90){\msmash{\tilde b}}
  \put(50,50){\circle{20}\msmash1}\put(40,65){\msmash{\tilde c}}
  \put(90,100){\vector(-1,0){80}}\put(50,110){\msmash{0}}
  \put(0,90){\vector(0,-1){80}}\put(10,50){\msmash{0,1}}
  \put(90,0){\vector(-1,0){80}}\put(50,10){\msmash{0}}
  \put(100,10){\vector(0,1){80}}\put(105,50){\msmash{1}}
  \put(43,43){\vector(-1,-1){36}}\put(25,35){\msmash{0}}
  \put(57,43){\vector(1,-1){36}}\put(75,35){\msmash{1}}
  \put(93,93){\vector(-1,-1){36}}\put(75,85){\msmash{1}}
  \put(-10,-10){\arc{20}{0}{4.712388981}}\put(0,-10){\vector(0,1){0}}\put(-10,-10){\msmash{0,1}}
\end{picture}
\caption{The Finite Transducers for $G$ and $\tilde G$}\label{fig:autoG}
\end{figure}

$G$ was originally defined in~\cite{grigorchuk:burnside-orig} as the
following dynamical system acting on the interval $[0,1]$ from which
rational dyadic points are removed:
\begin{align*}
  a(z)&=\begin{cases}z+\frac12&\text{ if }z<\frac12\\z-\frac12&\text{ if }z\ge\frac12,\end{cases}\\
  b(z)&=\begin{picture}(128,15)(0,-5)
    \put(0,0){\line(1,0){128}}
    \put(0,0){\line(0,1){10}}\put(32,7){\msmash a}
    \put(64,0){\line(0,1){10}}\put(80,7){\msmash a}
    \put(96,0){\line(0,1){10}}\put(104,7){\msmash 1}
    \put(112,0){\line(0,1){10}}\put(120,7){\msmash{a..}}
    \put(128,0){\line(0,1){10}}
    \put(0,-4){\tiny\msmash0}
    \put(64,-4){\tiny\msmash{\frac12}}
    \put(96,-4){\tiny\msmash{\frac34}}
    \put(112,-4){\tiny\msmash{\frac78}}
    \put(128,-4){\tiny\msmash1}
  \end{picture}\\
  c(z)&=\begin{picture}(128,10)
    \put(0,0){\line(1,0){128}}
    \put(0,0){\line(0,1){10}}\put(32,7){\msmash a}
    \put(64,0){\line(0,1){10}}\put(80,7){\msmash 1}
    \put(96,0){\line(0,1){10}}\put(104,7){\msmash a}
    \put(112,0){\line(0,1){10}}\put(120,7){\msmash{a..}}
    \put(128,0){\line(0,1){10}}
  \end{picture}\\
  d(z)&=\begin{picture}(128,10)
    \put(0,0){\line(1,0){128}}
    \put(0,0){\line(0,1){10}}\put(32,7){\msmash 1}
    \put(64,0){\line(0,1){10}}\put(80,7){\msmash a}
    \put(96,0){\line(0,1){10}}\put(104,7){\msmash a}
    \put(112,0){\line(0,1){10}}\put(120,7){\msmash{1..}}
    \put(128,0){\line(0,1){10}}
  \end{picture}
\end{align*}
Here the intervals represent $[0,1]$, with either $a$ or $1$ (the
identity transformation) acting on the described subintervals in a
similar way as $a$ or $1$ act on $[0,1]$. Finally, $G$ is also an
automata group, see the left graph in Figure~\ref{fig:autoG} (the
trivial and non-trivial elements of $\sym2$ are represented as
$1$ and $\epsilon$ and are used to label states).

Recall the map $\phi$ defined in~(\ref{eq:phi}); it restricts to an
embedding $\phi:H\to G\times G$ given by
\begin{equation}
  \phi:\begin{cases}b\to(a,c),\quad b^a\to(c,a)\\ c\to(a,d),\quad
    c^a\to(d,a)\\ d\to(1,b),\quad d^a\to(b,1),\end{cases}
  \label{eq:phiG}
\end{equation}
where $H=\st_G(1)=\langle b,c,d\rangle^G$ is an index-$2$ subgroup.
Consider also $K=\langle(ab)^2\rangle^G.$ Let $e$ be the infinite
sequence $0^\infty$; set $P_n=\st_G(0^n)$ and
$P=\st_G(e)=\cap_{n\ge0}P_n$ as above. Clearly $P_n$ has index $2^n$
in $G$, as $G$ acts transitively on $\Sigma^n$, and $P$ has infinite
index.

We note the following facts about $G$: it
\begin{itemize}
\item is an infinite torsion $2$-group;
\item is of intermediate growth, and therefore amenable (see
  Definition~\ref{def:amen});
\item is fractal, and regular branch on its subgroup $K$;
\item is just infinite;
\item has a recursive presentation with infinitely many relators;
  these relators are obtained as iterates of a substitution on a
  finite set of words~\cite{lysionok:pres};
\item is residually finite, and more precisely has a natural sequence
  of finite approximating quotients $G_n=G/\st_G(n)$, of order
  $2^{5\cdot2^{n-3}+2}$ for $n\ge3$ (and order $2^{2^n-1}$ for
  $n\le3$);
\item has a faithful action on the set $G/P$, of linear growth by
  Proposition~\ref{prop:polygrowth};
\item has a faithful action an $\partial\tree$ whose orbits are
  confinality classes.
\end{itemize}

The decomposition of $\rho_{G/P_n}$ in irreducibles is given by the
following lemma, combined with Lemma~\ref{lem:factorqr}:
\begin{lem}[\cite{bartholdi-g:parabolic}]
  $P_n$ has $n+1$ orbits in $\Sigma^n$; they are $0^n$ and the
  $0^i1\Sigma^{n-1-i}$ for $0\le i<n$. The orbits of $P$ in
  $\tree_\Sigma$ are the $0^i\Sigma^*$ for all $i\in\N$.
\end{lem}

\subsection{The Group $\tilde G$}
We describe briefly another fractal group, acting on the same tree
$\tree_2$ as $G$. More details appear in~\cite{bartholdi-g:parabolic}.
We denote again by $a$ the automorphism permuting the top two
branches, and let recursively $\tilde b$ be the automorphism acting as
$a$ on the right branch and $\tilde c$ on the left, $\tilde c$ be the
automorphism acting as $1$ on the right branch and $\tilde d$ on the
left, and $\tilde d$ be the automorphism acting as $1$ on the right
branch and $\tilde b$ on the left. In formul\ae,
\begin{alignat*}{2}
  \tilde b(0x\sigma)&=0\overline x\sigma,&\qquad\tilde b(1\sigma)&=1\tilde c(\sigma),\\
  \tilde c(0\sigma)&=0\sigma,&\qquad\tilde c(1\sigma)&=1\tilde d(\sigma),\\
  \tilde d(0\sigma)&=0\sigma,&\qquad\tilde d(1\sigma)&=1\tilde b(\sigma).
\end{alignat*}
Then $\tilde G$ is the group generated by $\{a,\tilde b,\tilde c,\tilde d\}$.
Clearly all these generators are of order $2$, and $\{\tilde b,\tilde
c,\tilde d\}$ is elementary abelian of order $8$. It can be defined
using the second automaton in Figure~\ref{fig:autoG}, or as the
dynamical system
\begin{align*}
  a(z)&=\begin{cases}z+\frac12&\text{ if }z<\frac12\\z-\frac12&\text{ if }z\ge\frac12,\end{cases}\\
  \tilde b(z)&=\begin{picture}(128,10)
    \put(0,0){\line(1,0){128}}
    \put(0,0){\line(0,1){10}}\put(32,7){\msmash a}
    \put(64,0){\line(0,1){10}}\put(80,7){\msmash 1}
    \put(96,0){\line(0,1){10}}\put(104,7){\msmash 1}
    \put(112,0){\line(0,1){10}}\put(120,7){\msmash{a..}}
    \put(128,0){\line(0,1){10}}
  \end{picture}\\
  \tilde c(z)&=\begin{picture}(128,10)
    \put(0,0){\line(1,0){128}}
    \put(0,0){\line(0,1){10}}\put(32,7){\msmash 1}
    \put(64,0){\line(0,1){10}}\put(80,7){\msmash 1}
    \put(96,0){\line(0,1){10}}\put(104,7){\msmash a}
    \put(112,0){\line(0,1){10}}\put(120,7){\msmash{1..}}
    \put(128,0){\line(0,1){10}}
  \end{picture}\\
  \tilde d(z)&=\begin{picture}(128,10)
    \put(0,0){\line(1,0){128}}
    \put(0,0){\line(0,1){10}}\put(32,7){\msmash 1}
    \put(64,0){\line(0,1){10}}\put(80,7){\msmash a}
    \put(96,0){\line(0,1){10}}\put(104,7){\msmash 1}
    \put(112,0){\line(0,1){10}}\put(120,7){\msmash{1..}}
    \put(128,0){\line(0,1){10}}
  \end{picture}
\end{align*}

Recall the map $\phi$ defined in~(\ref{eq:phi}); it restricts to an
embedding $\phi:\tilde H\to\tilde G\times\tilde G$ given by
\[\phi:\begin{cases}\tilde b\to(a,\tilde c),\quad\tilde b^a\to(\tilde c,a)\\ \tilde c\to(1,\tilde d),\quad\tilde c^a\to(d,1)\\ \tilde d\to(1,\tilde b),\quad\tilde d^a\to(\tilde b,1),\end{cases}\]
where $\tilde H=\st_{\tilde G}(1)=\langle\tilde b,\tilde c,\tilde
d\rangle^{\tilde G}$ is an index-$2$ subgroup. Consider also
$\tilde K=\langle (a\tilde b)^2,(a\tilde d)^2\rangle^{\tilde G}$.
Let again $e$ be the infinite sequence $0^\infty$; set $\tilde
P_n=\st_{\tilde G}(0^n)$ of index $2^n$, and $\tilde P=\st_{\tilde
  G}(e)=\cap_{n\ge0}\tilde P_n$ of infinite index.

We note the following facts about $G$, proved
in~\cite{bartholdi-g:parabolic}: it
\begin{itemize}
\item is an infinite group containing $G=\langle a,\tilde b\tilde
  c,\tilde c\tilde d,\tilde d\tilde b\rangle$ as an infinite-index
  subgroup, and has $2$-torsion elements as well as infinite-order 
  elements;
\item is of intermediate growth, and therefore is amenable;
\item is fractal, and regular branch on its subgroup $\tilde K$;
\item is just infinite;
\item has a recursive presentation with infinitely many relators;
  these relators are obtained as iterates of a substitution on a
  finite set of words. All of $\tilde G$'s relators have
  even length with respect to the generating set $\{a,\tilde b,\tilde
  c,\tilde d\}$, so its Cayley graph is bipartite;
\item is residually finite, and more precisely has a natural sequence
  of finite approximating quotients $\tilde G_n=\tilde G/\st_{\tilde
    G}(n)$, of order $2^{13\cdot2^{n-4}+2}$ for $n\ge 4$ (and order
  $2^{2^n-1}$ for $n\le4$).
\item has a faithful action on the set $\tilde G/\tilde P$, of
  linear growth by Proposition~\ref{prop:polygrowth};
\item has a faithful action an $\partial\tree$ whose orbits are
  confinality classes.
\end{itemize}

The decomposition of $\rho_{\tilde G/\tilde P_n}$ in irreducibles is
given by the following lemma, combined with Lemma~\ref{lem:factorqr}:
\begin{lem}
  $\tilde P_n$ has $n+1$ orbits in $\Sigma^n$; they are $0^n$ and the
  $0^i1\Sigma^{n-1-i}$ for $0\le i<n$. The orbits of $\tilde P$ in
  $\tree_\Sigma$ are the $0^i\Sigma^*$ for all $i\in\N$.
\end{lem}

\subsection{$GGS$ groups}\label{subs:GGS}
We next study three examples of a family of groups called $GGS$ groups
(the terminology was introduced by Gilbert
Baumslag~\cite{baumslag:cgt} and refers to Rostislav Grigorchuk,
Narain Gupta and Said Sidki). Let $p$ be a prime number. Denote $a$
the automorphism of $\tree_p$ permuting cyclically the top $p$
branches. Let
$\epsilon=(\epsilon_0,\dots,\epsilon_{p-2})\in(\Z/p)^{p-1}$. Define
recursively the automorphism $t_\epsilon$ of $\tree_p$ by
\[t(xy\sigma)=x(y+\epsilon_x)\sigma\text{ if }0\le x\le p-2,\qquad
t((p-1)\sigma)=(p-1)t(\sigma).\]
Then $G_\epsilon$ is the subgroup of $\aut(\tree_p)$ generated by
$\{a,t\}$.

The following results have their roots
in~\cite{grigorchuk:burnside-orig,gupta-s:infinitep} and are known as
part of folklore; for a proof see~\cite{grigorchuk:jibg}:
\begin{thm}
  $G_\epsilon$ is an infinite group if an only if
  $\epsilon\neq(0,\dots,0)$. It is a torsion group if and only if
  $\sum\epsilon_i=0$.
\end{thm}

\begin{thm}[\cite{bartholdi:ggs}]
  For all $\epsilon$, the group $G_\epsilon$ is of subexponential
  growth, and therefore is amenable.
\end{thm}

Let $G=G_\epsilon$ be as above. Let $e=(p-1)^\infty$ be the rightmost
path in $\tree_p$ and let $P=\st_G(e)$.
Proposition~\ref{prop:polygrowth} applies, so $G/P$ has polynomial
growth of degree at most $\log_2(p)$.

\subsection{The Group $\Gamma$}
\begin{figure}
\begin{picture}(150,150)(-25,-25)
  \put(0,0){\circle{20}\msmash\varepsilon}\put(-10,15){\msmash a}
  \put(100,0){\circle{20}\msmash1}
  \put(110,-10){\arc{20}{4.712388981}{9.424777962}}\put(110,0){\vector(-1,0){0}}\put(111,-11){\msmash\Sigma}
  \put(50,100){\circle{20}\msmash1}\put(38,92){\msmash s}
  \put(10,0){\vector(1,0){80}}\put(50,-5){\msmash\Sigma}
  \put(45.5,91){\vector(-1,-2){41}}\put(19,50){\msmash0}
  \put(54.5,91){\vector(1,-2){41}}\put(81,50){\msmash1} 
  \put(50,114){\arc{20}{2.356194491}{7.068583472}}\put(43,107){\vector(1,-1){0}}
  \put(50,118){\msmash2}
\end{picture}
\qquad
\begin{picture}(50,150)(-25,-25)
  \put(0,50){\circle{20}\msmash\varepsilon}\put(-10,40){\msmash a}
  \put(0,100){\circle{20}\msmash1}\put(-10,90){\msmash t}
  \put(0,0){\circle{20}\msmash1}
  \put(0,90){\vector(0,-1){30}}\put(10,75){\msmash{0,1}}
  \put(0,40){\vector(0,-1){30}}\put(5,25){\msmash\Sigma}
  \put(0,-14){\arc{20}{5.497787145}{10.21017613}}\put(7,-7){\vector(-1,1){0}}
  \put(0,-16){\msmash\Sigma}
  \put(0,114){\arc{20}{2.356194491}{7.068583472}}\put(-7,107){\vector(1,-1){0}}
  \put(0,118){\msmash2}
\end{picture}
\qquad
\begin{picture}(150,150)(-25,-25)
  \put(0,0){\circle{20}\msmash\varepsilon}\put(-10,15){\msmash a}
  \put(0,100){\circle{20}\msmash1}\put(-10,90){\msmash r}
  \put(100,0){\circle{20}\msmash1}
  \put(100,100){\circle{20}\msmash{\varepsilon^2}}
  \put(0,90){\vector(0,-1){80}}\put(-5,50){\msmash0}
  \put(100,90){\vector(0,-1){80}}\put(105,50){\msmash\Sigma}
  \put(10,0){\vector(1,0){80}}\put(50,-5){\msmash\Sigma}
  \put(10,100){\vector(1,0){80}}\put(50,107){\msmash1}
  \put(-10,110){\arc{20}{1.570796327}{6.283185308}}\put(-10,100){\vector(1,0){0}}\put(-11,111){\msmash2}
  \put(110,-10){\arc{20}{4.712388981}{9.424777962}}\put(110,0){\vector(-1,0){0}}\put(111,-11){\msmash\Sigma}
\end{picture}
\caption{The Automata for $\Gamma$, $\overline\Gamma$ and $\doverline\Gamma$}\label{fig:autoGamma}
\end{figure}
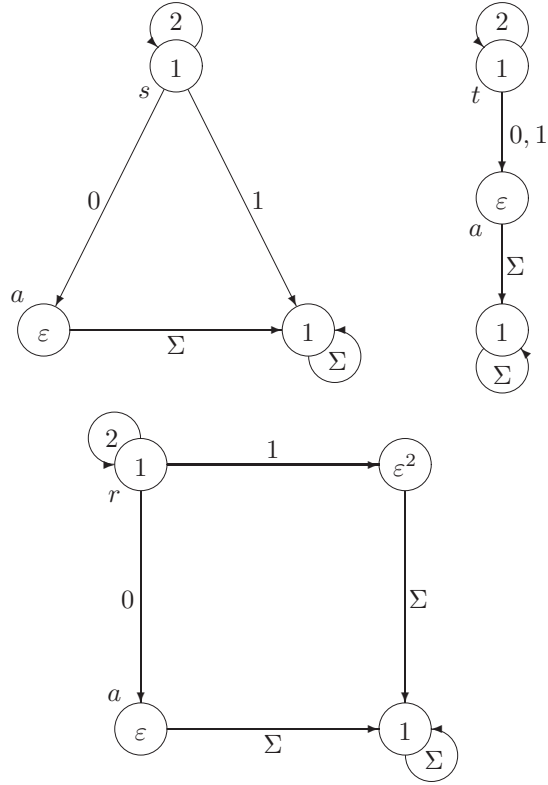

As always denote by $a$ the rooted automorphism of $\tree_3$ permuting
cyclically the top three branches. Let $s$ be the automorphism of
$\tree_3$ defined recursively by
\[s(0x\sigma)=0\overline x\sigma,\qquad s(1x\sigma)=1x\sigma,\qquad
s(2\sigma)=2s(\sigma).\]
Then $\Gamma$ is the subgroup of
$\aut(\tree_3)$ generated by $\{a,s\}$; its growth was studied by
Jacek Fabrykowski and Narain Gupta~\cite{gupta-f:growth2}. It can also
be defined using automata, as in Figure~\ref{fig:autoGamma}, or as the
dynamical system
\begin{align*}
  a(z)&=\begin{cases}z+\frac13&\text{ if }z<\frac23\\z-\frac23&\text{ if }z\ge\frac23,\end{cases}\\
  s(z)&=\begin{picture}(144,15)(0,-5)
    \put(0,0){\line(1,0){144}}
    \put(0,0){\line(0,1){10}}\put(24,7){\msmash a}
    \put(48,0){\line(0,1){10}}\put(72,7){\msmash 1}
    \put(96,0){\line(0,1){10}}\put(104,7){\msmash a}
    \put(112,0){\line(0,1){10}}\put(120,7){\msmash 1}
    \put(128,0){\line(0,1){10}}\put(136,7){\msmash{a..}}
    \put(144,0){\line(0,1){10}}
    \put(0,-4){\tiny\msmash0}
    \put(48,-4){\tiny\msmash{\frac13}}
    \put(96,-4){\tiny\msmash{\frac23}}
    \put(112,-4){\tiny\msmash{\frac79}}
    \put(128,-4){\tiny\msmash{\frac89}}
    \put(144,-4){\tiny\msmash1}
  \end{picture}\\
\end{align*}

Set $H=\st_\Gamma(1)=\langle s\rangle^\Gamma$; then
$\phi:H\to\Gamma\times\Gamma\times\Gamma$ can be expressed as
\begin{equation}
  \phi:\begin{cases}s\to(a,1,s),\quad s^a\to(s,a,1),\quad s^{a^2}\to(1,s,a).\end{cases}
  \label{eq:phiGamma}
\end{equation}

Define the elements $x=as$, $y=sa$ of $\Gamma$, and let $K$ be the
subgroup of $\Gamma$ generated by $x$ and $y$. Then $K$ is normal in
$\Gamma$, because $y^s=x^{-1}y^{-1}$, $y^{a^{-1}}=y^{-1}x^{-1}$,
$y^{s^{-1}}=y^a=a$, and similar relations hold for conjugates of $x$.
Moreover $K$ is of index $3$ in $\Gamma$, with transversal $\langle
a\rangle$. Let $L$ be the subgroup of $K$ generated by $[K,K]$ and cubes
in $K$.

\begin{thm}[\cite{bartholdi-g:parabolic}]
  $\Gamma$ is a regular branch, fractal group.  The subgroup $K$ of $\Gamma$
  is torsion-free; thus $\Gamma$ is virtually torsion-free.  The
  finite quotients $\Gamma_n=\Gamma/H_n$ of $\Gamma$ have order
  $3^{3^{n-1}+1}$ for $n\ge2$, and $3$ for $n=1$.
\end{thm}

\subsection{The Group $\overline\Gamma$}
Let $a$ be as in the previous subsection, and let $t$ be the
automorphism of $\tree_3$ defined recursively by
\[t(0x\sigma)=0\overline x\sigma,\qquad t(1x\sigma)=1\overline
x\sigma,\qquad t(2\sigma)=2t(\sigma).\]
Then $\overline\Gamma$ is the subgroup of $\aut(\tree_3)$ generated by
$\{a,t\}$. The associated dynamical system is
\begin{align*}
  a(z)&=\begin{cases}z+\frac13&\text{ if }z<\frac23\\z-\frac23&\text{ if }z\ge\frac23,\end{cases}\\
  t(z)&=\begin{picture}(144,10)
    \put(0,0){\line(1,0){144}}
    \put(0,0){\line(0,1){10}}\put(24,7){\msmash a}
    \put(48,0){\line(0,1){10}}\put(72,7){\msmash a}
    \put(96,0){\line(0,1){10}}\put(104,7){\msmash a}
    \put(112,0){\line(0,1){10}}\put(120,7){\msmash a}
    \put(128,0){\line(0,1){10}}\put(136,7){\msmash{a..}}
    \put(144,0){\line(0,1){10}}
  \end{picture}.\\
\end{align*}

Set $H=\st_{\overline\Gamma}(1)=\langle t\rangle^{\overline\Gamma}$; then
$\phi:H\to\overline\Gamma\times\overline\Gamma\times\overline\Gamma$
can be expressed as
\[\phi:\begin{cases}t\to(a,a,t),\quad t^a\to(t,a,a),\quad t^{a^2}\to(a,t,a).\end{cases}\]

Define the elements $x=ta^{-1}$, $y=a^{-1}t$ of $\overline\Gamma$, and
let $\overline K$ be the subgroup of $\overline\Gamma$ generated by
$x$ and $y$.  Then $\overline K$ is normal in $\overline\Gamma$,
because $x^t=y^{-1}x^{-1}$, $x^a=x^{-1}y^{-1}$,
$x^{t^{-1}}=x^{a^{-1}}=y$, and similar relations hold for conjugates
of $y$.  Moreover $\overline K$ is of index $3$ in $\overline\Gamma$,
with transversal $\langle a\rangle$.

\begin{thm}[\cite{bartholdi-g:parabolic}]\label{prop:Gammafractal}
  $\overline\Gamma$ is a fractal group and is weak branch, but not
  branch.  The subgroup $\overline K$ of $\overline\Gamma$ is
  torsion-free; thus $\overline\Gamma$ is virtually torsion-free.  The
  finite quotients $\overline\Gamma_n=\overline\Gamma/H_n$ of
  $\overline\Gamma$ have order $3^{\frac14(3^n+2n+3)}$ for $n\ge2$,
  and $3^{\frac12(3^n-1)}$ for $n\le2$.
\end{thm}

\subsection{The Group $\doverline\Gamma$}
Let again $a$ denote the automorphism of $\tree_3$ permuting cyclically
the top three branches. Let now $r$ be the automorphism of $\tree_3$
defined recursively by
\[r(0x\sigma)=0\overline x\sigma,\qquad r(1x\sigma)=1\doverline x\sigma,\qquad
r(2\sigma)=2r(\sigma).\]

Then $\doverline\Gamma$ is the subgroup of $\aut(\tree_3)$
generated by $\{a,r\}$; it was studied by Gupta and Sidki~\cite{gupta-s:3group,gupta-s:infinitep,sidki:subgroups,sidki:presentation}. The associated dynamical system is
\begin{align*}
  a(z)&=\begin{cases}z+\frac13&\text{ if }z<\frac23\\z-\frac23&\text{ if }z\ge\frac23,\end{cases}\\
  r(z)&=\begin{picture}(144,10)
    \put(0,0){\line(1,0){144}}
    \put(0,0){\line(0,1){10}}\put(24,7){\msmash a}
    \put(48,0){\line(0,1){10}}\put(72,7){\msmash{a^2}}
    \put(96,0){\line(0,1){10}}\put(104,7){\msmash a}
    \put(112,0){\line(0,1){10}}\put(120,7){\msmash{a^2}}
    \put(128,0){\line(0,1){10}}\put(136,7){\msmash{a..}}
    \put(144,0){\line(0,1){10}}
  \end{picture}.\\
\end{align*}

Set $H=\st_{\doverline\Gamma}(1)=\langle r\rangle^{\doverline\Gamma}$; then
$\phi:H\to\doverline\Gamma\times\doverline\Gamma\times\doverline\Gamma$
con be expressed as
\[\phi:\begin{cases}r\to(a,a^2,r),\quad r^a\to(r,a,a^2),\quad r^{a^2}\to(a^2,r,a).\end{cases}\]

\begin{thm}[\cite{bartholdi-g:parabolic}]
  $\doverline\Gamma$ is a just infinite torsion $3$-group. It is branch
  and fractal. The finite quotients
  $\doverline\Gamma_n=\doverline\Gamma/H_n$ of
  $\Gamma$ have order $3^{2\cdot3^{n-2}+1}$ for $n\ge2$, and $3$ for
  $n=1$.
\end{thm}

\section{Unitary Representations and Hecke type Operators}\label{sec:urho}
The five groups introduced in the previous section share the property
of acting faithfully on a regular rooted tree; natural representations
arise from this fact. In this chapter we suppose $G$ is any group
acting level-transitively on a regular tree.

We defined in Subsection~\ref{subs:ds} the boundary $\partial\tree$ of
the tree on which $G$ acts. Since $G$ preserves the uniform measure on
this boundary, we have a unitary representation $\pi$ of $G$ in
$L^2(\partial\tree,\nu)$, or equivalently in $L^2([0,1],m)$. Let
$\hilb_n$ be the subspace of $L^2(\partial\tree,\nu)$ spanned by the
characteristic functions $\chi_\sigma$ of the rays $e$ starting by
$\sigma$, for all $\sigma\in\Sigma^n$. It is of dimension $d^n$, and
can equivalently be seen as spanned by the characteristic functions in
$L^2([0,1],m)$ of intervals of the form $[(i-1)d^{-n},id^{-n}]$, $1\le
i\le d^n$. These $\hilb_n$ are invariant subspaces, and afford
representations $\pi_n=\pi_{|\hilb_n}$. As clearly $\pi_{n-1}$ is a
subrepresentation of $\pi_n$, we set $\pi^\perp_n=\pi_n\ominus\pi_{n-1}$,
so that $\pi=\oplus_{n=0}^\infty\pi^\perp_n$.

Denote by $\rho_{G/H}$ the quasi-regular representation of $G$ in
$\ell^2(G/H)$ and by $\rho_{G/H_n}$ the finite-dimensional
representations of $G$ in $\ell^2(G/H_n)$. Since $G$ is
level-transitive, the representations $\pi_n$ and $\rho_{G/H_n}$ are
unitary equivalent.

\begin{defn}\label{defn:sch}
  Let $G$ be a group generated by a set $S$ and $H$ a subgroup of
  $G$. The \emdef{Schreier graph} $\sch(G,H,S)$ of $G/H$ is the
  directed graph on the edge set $G/H$, with for every $s\in S$ and
  every $gH\in G/H$ an edge from $gH$ to $sgH$. The base point of
  $\sch(G,H,S)$ is the coset $H$.
\end{defn}
Note that $\sch(G,1,S)$ is the Cayley graph of $G$ relative to $S$. It
may happen that $\sch(G,P,S)$ have loops and multiple edges even if
$S$ is disjoint from $H$. Schreier graphs are $|S|$-regular graphs,
and any degree-regular graph $\gf$ containing a $1$-factor (i.e.\ a
regular subgraph of degree $1$; there is always one if $\gf$ has even
degree) is a Schreier graph~\cite[Theorem~5.4]{lubotzky:cayley}.

\begin{defn}
  Let $G$ be a group generated by a finite symmetric set $S$. The
  \emdef{spectrum} $\spec(\tau)$ of a representation $\tau:G\to\uhilb$
  with respect to the given set of generators is the spectrum of
  $\Delta_\tau=\sum_{s\in S}\tau(s)$ seen as an bounded operator on
  $\hilb$.
\end{defn}
As the vertices of $\sch(G,H,S)$ coincide with the set $G/H$, it is
easy to see that $\Delta_\tau/|S|-1$ is unitary equivalent to the
Laplacian operator on $\sch(G,H,S)$, and therefore has same spectrum.

\begin{defn}\label{def:amen}
  Let $G$ be a group acting on a set $X$. This action is
  \emdef{amenable} in the sense of von
  Neumann~\cite{vneumann:masses} if there exists a finitely additive
  measure $\mu$ on $X$, invariant under the action of $G$, with
  $\mu(X)=1$.

  A group $G$ is \emdef{amenable} if its action on itself by
  left-multiplication is amenable.
\end{defn}

Amenability can be tested using the following criterion, due to F\o
lner for the regular action~\cite{folner:banach} (see
also~\cite{ceccherini-:amen} and the literature cited there):
\begin{thm}
  Assume the group $G$ acts on a discrete set $X$. Then the action is
  amenable if and only for every if for every $\lambda>0$ and every
  $g\in G$ there exists a finite subset $F\subset X$ such that $|F\triangle
  gF|<\lambda|F|$, where $\triangle$ denotes symmetric difference and
  $|\cdot|$ cardinality.
\end{thm}
Using this criterion it is easy to see that any $G$-space $X$ of
subexponential growth is amenable. In particular the $G$-spaces $G/P$
are of polynomial growth when the conditions of
Proposition~\ref{prop:polygrowth} are fulfilled, and therefore are
amenable. The following result belongs to the common lore, and we
don't know of a reference to its proof. Our attention was drawn to it
by Marc Burger and Alain Valette:
\begin{prop}\label{prop:burger}
  Let $H<G$ be any subgroup. Then the quasi-regular representation
  $\rho_{G/H}$ is weakly contained in $\rho_G$ if and only if
  $H$ is amenable.
\end{prop}
\begin{proof}
  If $H$ is amenable, then the trivial one-dimensional representation
  $1_H$ of $H$ is weakly contained in $\rho_H$. Inducing up,
  $\rho_{G/H}=\operatorname{Ind}_H^G1_H$ is weakly contained in
  $\rho_G=\operatorname{Ind}_H^G(\rho_H)$.

  Conversely, if $\rho_{G/H}$ is weakly contained in $\rho_G$, we
  obtain by restricting to $H$ that ${\rho_{G/H}}_{|H}$ is weakly
  contained in ${\rho_G}_{|H}$. Now $1_H$ is a subrepresentation of
  ${\rho_{G/H}}_{|H}$, because the Dirac mass at $H$ is $H$-fixed; and
  ${\rho_G}_{|H}=[G:H]\rho_H$ by Frobenius reciprocity. It follows
  that $1_H$ is weakly contained in $[G:H]\rho_H$, and therefore that
  $H$ is amenable.
\end{proof}

The following statements allows one to compare spectra of diverse representations.
\begin{thm}
  Let $G$ be a group acting on a regular rooted tree, and let $\pi$,
  $\pi_n$ and~$\pi^\perp_n$ be as above.
  \begin{enumerate}
  \item If $G$ is weak branch, then $\rho_{G/P}$ is an irreducible
    representation of infinite dimension.
  \item $\pi$ is a reducible representation of infinite dimension
    whose irreducible components are precisely those of the $\pi^\perp_n$
    (and thus are all finite-dimensional). Moreover
    \[\spec(\pi)=\overline{\bigcup_{n\ge0}\spec(\pi_n)}=\overline{\bigcup_{n\ge0}\spec(\pi^\perp_n)}.\]
  \item The spectrum of $\rho_{G/P}$ is contained in
    $\overline{\cup_{n\ge0}\spec(\rho_{G/P_n})}=\overline{\cup_{n\ge0}\spec(\pi_n)}$,
    and thus is contained in the spectrum of $\pi$. If moreover either
    $P$ or $G/P$ are amenable, these spectra coincide, and if $P$ is
    amenable, they are contained in the spectrum of $\rho_G$:
    \[\spec(\rho_{G/P})=\spec(\pi)=\overline{\bigcup_{n\ge0}\spec(\pi_n)}\subseteq\spec(\rho_G).\]
  \item $\Delta_\pi$ has a pure-point spectrum, and its spectral radius
    $r(\Delta_\pi)=s\in\R$ is an eigenvalue, while the spectral
    radius $r(\Delta_{\rho_{G/P}})$ is not an eigenvalue of
    $\Delta_{\rho_{G/P}}$. Thus $\Delta_{\rho_{G/P}}$ and $\Delta_\pi$
    are different operators having the same spectrum.
  \end{enumerate}
\end{thm}
\begin{proof}
  The first statement follows from Mackey's theorem~\ref{thm:mackey}
  and Proposition~\ref{prop:comm}. The second holds because $\pi$
  splits as the direct sum of the $\pi^\perp_n$'s, and is weakly
  equivalent to the direct sum of the $\pi_n$'s.
  
  It was mentioned that $\rho_{G/P_n}$ and $\pi_n$ are equivalent:
  they are both finite-dimensional and act on $G$-equivalent sets,
  namely $G/H_n$ and $\Sigma^n$. The third statement then follows from
  Proposition~\ref{prop:burger} and Propositions~\ref{prop:spinc}
  and~\ref{prop:speq} below.

  It is obvious that $s$ is an eigenvalue of $\Delta_\pi$ with
  constant eigenfunction. Now the fourth follows from
  Proposition~5 in~\cite{grigorchuk-z:infinite}.
\end{proof}

Since all of our example groups are amenable, the spectra computed in
the next section are included in $\spec(\rho_G)$. Moreover, as $\tilde
G$ has a bipartite Cayley graph, $\spec(\rho_{\tilde G})$ is symmetrical
about $1$ and contains $[0,4]$ (as we shall show in
Section~\ref{subs:spgt}), so is $[-4,4]$.

We finish this subsection by turning to a question of Mark
Ka\'c~\cite{kac:drum}: ``Can one hear the shape of a drum?'' This
question was answered in the negative in~\cite{gordon-w-w:drum}, and
we here answer by the negative to a related question: ``Can one hear a
representation?'' Indeed $\rho_{G/P}$ and $\pi$ have same spectrum
(i.e.\ cannot be distinguished by hearing), but are not equivalent.
Furthermore, if $G$ is a branched group, there are uncountably many
nonequivalent representations within
$\{\rho_{G/\st_G(e)}|\,e\in\partial\tree\}$, as is shown
in~\cite{bartholdi-g:parabolic}.

The same question may be asked for graphs: ``are there two
non-isomorphic graphs with same spectrum?'' There are finite examples,
obtained through the notion of \emph{Sunada
  pair}~\cite{lubotzky:cayley}. C\'edric B\'eguin, Alain Valette and
Andrzej \.Zuk produced the following example in~\cite{beguin-:spectrum}:
let $\Gamma$ be the integer Heisenberg group (free $2$-step nilpotent
on $2$ generators $x,y$). Then $\Delta=x+x^{-1}+y+y^{-1}$ has spectrum
$[-2,2]$, which is also the spectrum of $\Z^2$ for an independent
generating set. As a consequence, their Cayley graphs have same
spectrum, but are not quasi-isometric (they do not have the same growth).

Using the result of Nigel Higson and Gennadi
Kasparov~\cite{higson-k:bc} (giving a partial positive answer to the
Baum-Connes conjecture), we may infer the following
\begin{prop}
  Let $\Gamma$ be a torsion-free amenable group with finite generating
  set $S=S^{-1}$ such that there is a map $\phi:\Gamma\to\Z/2\Z$ with
  $\phi(S)=\{1\}$.  Then
  \[\spec(\sum_{s\in S}\rho(s))=[-|S|,|S|].\]
\end{prop}
In particular, there are countably many non-quasi-isometric graphs
with the same spectrum, including the graphs of $\Z^d$, of free
nilpotent groups and of suitable torsion-free groups of intermediate
growth (for the first examples, see~\cite{grigorchuk:pgps}).
\begin{proof}
  Since $\Gamma$ is amenable, $|S|$ is in its spectrum. By the
  existence of $\phi$, the Cayley graph of $\Gamma$ with respect to
  $S$ is two-colourable (i.e.\ bipartite), so its spectrum is
  symmetrical, and therefore contains $-|S|$. By the Baum-Connes
  conjecture, proved for this case by Higson and Kasparov, the
  $C^*$-algebra $C^*_r(\rho(\Gamma))$ contains no idempotents. It
  follows by functional integration that $\rho$'s spectrum is
  connected, so is $[-|S|,|S|]$ as claimed.
\end{proof}


\subsection{Approximations of Operators and Spectra}\label{subs:approx}
We now prove the claimed inclusions of spectra, in part relying on
$C^*$-algebraic results.
\begin{prop}\label{prop:spinc}
  Let $\{H_n\}_{n\ge0}$ be a descending family of finite-index
  subgroups of $G$, and set $H=\bigcap_{n\ge0}H_n$. Let $\tau$ and
  $\tau_n$ be the quasi-regular representations of $G$ on $G/H$ and
  $G/H_n$ respectively. Then
  \[\spec\tau\subseteq\overline{\bigcup_{n\ge0}\spec \tau_n}.\]
\end{prop}
\begin{proof}
  Let $\sch$ and $\sch_n$ be the Schreier graphs $\sch(G,H,S)$ and
  $\sch(G,H_n,S)$ respectively, and mark in them the vertices $H$ and
  $H_n$. Then $\sch_n\underset{n\to\infty}\Longrightarrow\sch$ in the
  sense of~\cite{grigorchuk-z:infinite}, that is, coincide with $\sch$
  in ever increasing balls centered at $H_n$; therefore, for their
  associated spectral measures,
  $\sigma_n(\lambda)\underset{n\to\infty}\Longrightarrow\sigma(\lambda)$,
  in the sense of weak convergence. Therefore the support of $d\sigma$
  is contained in the closure of the union of the supports of
  $d\sigma_n$, and the proposition follows.
\end{proof}

\begin{proof}[Alternate Proof]
  The quasi-regular representation $\tau$ is weakly contained in
  $\oplus_{n\ge0}\tau_n$, by ``approximation of
  coefficients''~\cite[Theorem~3.4.9]{dixmier:c*algebras}. Indeed, as
  $\delta_H$, the Dirac function at $H$, is a cyclic vector for
  $\tau$, it is enough to see that the function
  \[g\mapsto\langle\tau(g)\delta_H|\delta_H\rangle=\begin{cases}1&\text{if }g\in H\\0&\text{otherwise}\end{cases}\]
  can be pointwise approximated on $G$ be coefficients of the
  $\tau_n$'s. But
  \[\langle\tau_n(g)\delta_{H_n}|\delta_{H_n}\rangle=\begin{cases}1&\text{if }g\in H_n\\0&\text{otherwise}\end{cases},\]
  so
  \[\lim_{n\to\infty}\langle\tau_n(g)\delta_{H_n}|\delta_{H_n}\rangle=\langle\tau(g)\delta_H|\delta_H\rangle.\]
  We then have a surjection
  $C^*(\oplus_{n\ge0}\tau_n)\twoheadrightarrow C^*(\tau)$, where
  $C^*(\rho)$ is the $C^*$-algebra generated by the image of
  $\rho$. The spectrum inclusions follow.
\end{proof}

\begin{prop}\label{prop:speq}
  Let $\{H_n\}_{n\ge0}$ be a descending family of finite-index
  subgroups of $G$, and set $H=\bigcup_{n\ge0}H_n$. Let $\tau$ and
  $\tau_n$ be the quasi-regular representations of $G$ on $G/H$ and
  $G/H_n$ respectively. Assume moreover that the action of $G$ on
  $G/H$ is amenable, i.e.\ that $\sch(G,H,S)$ is amenable. Then
  \[\spec\tau=\overline{\bigcup_{n\ge0}\spec \tau_n}.\]
\end{prop}
\begin{proof}
  Let $\sch$ be the Schreier graph $\sch(G,H,S)$. Choose
  $\mu\in\spec\tau_n$, with a corresponding eigenvector $v:G/H_n\to\C$.
  As $\sch$ is amenable, there is a F\o lner sequence
  $\{F_k\}$ in $\sch$, i.e.\ a family of subsets of $V(\sch)$ with
  $|F_k|/|\partial F_k|\to0$ as $k\to\infty$.
  Define now functions $v_k$ on $\sch$,
  \[v_k(gH)=\begin{cases}v(gH_n) &\text{ if }gH\in F_k,\\0&\text{ otherwise.}\end{cases}\]
  Then $\|\Delta_\tau v_k-\mu v_k\|\le|\partial F_k|\cdot\|v\|$. If we let
  $\Omega$ be a fundamental domain of $G/H_n$ in $\sch$, of diameter
  $\delta$, and let $\partial_\delta(F_k)$ denote the $\delta$-neighbourhood
  of $F_k$, then we have $|\Omega|\cdot\|v_k\|\ge(|F_k|-|\partial_\delta
  F_k|)\|v\|$, whence
  \[\left\|\Delta_\tau\frac{v_k}{\|v_k\|}-\mu\frac{v_k}{\|v_k\|}\right\|\le
  \frac{|\Omega|\cdot|\partial F_k|}{|F_k|-|\partial_\delta F_k|}\to0,\]
  so $\mu\in\spec\tau$.
\end{proof}

%
%

\subsection{Spectral Measures}\label{subs:sm}
In parallel to the computation of spectra, interesting questions arise
in relation to properties of spectral measures associated to Markovian
operators. There are two approaches to spectral measures, one via a
solution to the moment problem and one via trace states in von Neumann
algebras.

Let $\gf$ be a graph, and let $M$ be its Markovian operator. For
any $x,y\in V(\gf)$ and $n\in\N$ let $p_{x,y}^n$ be the probability
that a simple random walk starting at $x$ be at $y$ after $n$ steps. Recall
that if $M_v$ be the characteristic vector of the vertex $v$, then
\[p_{x,y}^n = \langle M^n\delta_x|\delta_y\rangle.\]
Define the spectral measures $\sigma_{x,y}$ by
\[p_{x,y}^n = \int_{-1}^1\lambda^nd\sigma_{x,y}(\lambda)\qquad \forall
n\in\N,\]
or equivalently as
\[\sigma_{x,y}(\lambda) = \langle M(\lambda)\delta_x|\delta_y\rangle,\]
where $M(\lambda)$ is the spectral decomposition of $M$ (the operator
coinciding with $M$ on eigenfunctions whose eigenvalue is at most
$\lambda$). Set $\sigma_x = \sigma_{x,x}$. These
measures are called the \emdef{Kesten spectral measures}, as they were
introduced in~\cite{kesten:rwalks}.
\begin{prop}
  If $\gf$ is connected, then all the measures $\sigma_{x,y}$ are
  equivalent.
\end{prop}
Therefore there is only one type of Kesten spectral measure, up to
equivalence.

Now let $N$ be a von Neumann algebra with a finite state $\tau$. For
any self-adjoint element $a\in N$ one can define a spectral measure
$\tau_a$ by $\tau_a(\lambda)=\tau(a(\lambda))$, where $a(\lambda)$ is
the spectral decomposition of $a$. We shall call $\tau_a$ the
\emph{von Neumann spectral measure} associated to $a$.

For instance, such a situation appears if $N$ is finite-dimensional,
or more generally of finite type. Another important example is when
$N$ is the von Neumann algebra generated by the left-regular
representation of a group. In this case the von Neumann spectral
measure is given by $\tau_a(\lambda)=\langle
a(\lambda)\delta_1,\delta_1\rangle$, where $\delta_1$ is the Dirac
function in the identity of the group. Therefore in this case the
von Neumann and Kesten spectral measures of a Markovian operator
coincide. A similar situation occurs if $N$ is a hyperfinite algebra
of type $II_1$, but there $\delta_1$ should be replaced by any cyclic
vector. If $N$ is any algebra of type $II_1$ there is a canonical
trace of the form $\tau(a)=\sum_i\langle ax_i,x_i\rangle$, for some
(in general infinite) sequence of vectors $(x_i)$.

Let $P$ be a subgroup of a group $G$ and suppose that the von Neumann
algebra generated by $\rho_{G/P}$ has finite trace. Then for any
symmetric system $S$ of generators of $G$ the von Neumann spectral
measure $\tau_a=\tau_{\sum_{s\in S}s}$ can be defined. We also call
$\tau_a$ the von Neumann spectral measure of the Schreier graph
$\gf=\sch(G,P,S)$.

In case $\gf$ is finite, the von Neumann spectral measure is just a
``histogram'', counting in any given interval the ``average number of
eigenvalues'' that belong to it.

Suppose now that $P=\cap_{n=1}^\infty P_n$, where the $P_n$ are
subgroups of finite index of $G$.  Then the sequence of finite graphs
$\gf_n=\sch(G,P_n,S)$ converges to the graph $\gf$ in the sense
of~\cite{grigorchuk-z:infinite}, and the following statement holds:
\begin{prop}[\cite{grigorchuk-z:infinite}]
  Let $\sigma_n$ and $\sigma$ be the Kesten spectral measures of the
  Schreier graphs $\gf_n$ and $\gf$ based at $P_n$ and $P$
  respectively.

  Then we have $\sigma_n(\lambda)\Longrightarrow\sigma(\lambda)$ in
  the sense of weak convergence.
\end{prop}

In case the subgroups $P_n$ are normal and $P=1$, we deduce from this
proposition the convergence of the corresponding spectral measures
$\tau_n$ to $\tau$, since in this case $\sigma_n=\tau_n$ and
$\sigma=\tau$. These results were obtained by Wolfgang
L\"uck~\cite{luck:l2} and Michael Farber~\cite{farber:geometry} using
different methods.

It is interesting to study the conditions under which the spectral
measures $\tau_n$ of finite graphs $\sch(G,P_n,S)$ converge to some
limit $\tau_*$ (which we call the \emph{empiric} spectral measure),
and, in case the von Neumann spectral measure $\tau$ is well defined
for a graph $\sch(G,P,S)$ (for instance, for a quasi-regular
representation $\rho_{G/P}$ generating a von Neumann algebra of finite
type), under which conditions we have $\tau=\tau_*$. Also, in the last
case, when do we have $\sigma=\tau$?

Unfortunately, in our situation, the von Neumann algebra generated by
$\rho_{G/P}$ is the algebra of all bounded operators so it has no a
good state.

More investigations should be done in order to clarify the meaning of
our computations of empiric spectral measures in the examples that follow.

\subsection{Operator Recursions}\label{subs:rec}
Let $\hilb$ be an infinite-dimensional Hilbert space, and suppose
$\Phi:\hilb\to\hilb\oplus\dots\oplus\hilb$ is an isomorphism, where
the domain of $\Phi$ is a sum of $d\ge2$ copies of $\hilb$. Let $S$ be
a finite subset of $\uhilb$, and suppose that for all $s\in S$, if we
write $\Phi^{-1}s\Phi$ as an operator matrix
$(s_{i,j})_{i,j\in\{1,\dots,d\}}$ where the $s_{i,j}$ are operators in
$\hilb$, then $s_{i,j}\in S\cup\{0,1\}$.

This is precisely the setting in which we will compute the spectra of
our five example groups: for $G$, we have $d=2$ and $S=\{a,b,c,d\}$ with
\begin{gather*}
  a = \begin{pmatrix}0&1\\1&0\end{pmatrix},\qquad b = \begin{pmatrix}a&0\\0&c\end{pmatrix},\\
  c = \begin{pmatrix}a&0\\0&d\end{pmatrix},\qquad d = \begin{pmatrix}1&0\\0&b\end{pmatrix}.
\end{gather*}
For $\tilde G$, we also have $d=2$, and $S=\{a,\tilde b,\tilde
c,\tilde d\}$ given by
\begin{gather*}
  b = \begin{pmatrix}a&0\\0&c\end{pmatrix},\qquad c = \begin{pmatrix}1&0\\0&d\end{pmatrix},\qquad d = \begin{pmatrix}1&0\\0&b\end{pmatrix}.
\end{gather*}
For $\Gamma=\langle a,s\rangle$, $\overline\Gamma=\langle a,t\rangle$
and $\doverline\Gamma=\langle a,r\rangle$, we have $d=3$ and
\begin{gather*}
  a = \begin{pmatrix}0&1&0\\0&0&1\\1&0&0\end{pmatrix},\qquad
  s = \begin{pmatrix}a&0&0\\0&1&0\\0&0&s\end{pmatrix},\\
  t = \begin{pmatrix}a&0&0\\0&a&0\\0&0&t\end{pmatrix},\qquad
  r = \begin{pmatrix}a&0&0\\0&a^2&0\\0&0&r\end{pmatrix}.
\end{gather*}

Each of these operators is unitary. The families $S=\{a,b,c,d\},\dots$
generate subgroups $G(S)$ of $\uhilb$.  The choice of the isomorphism
$\Phi$ defines a unitary representation of $\langle S\rangle$.

We note, however, that the expression of each operator as a matrix of
operators does not uniquely determine the operator, in the sense that
different isomorphisms $\Phi$ can yield non-isomorphic operators
satisfying the same recursions. Even if $\Phi$ is fixed, it may happen
that different operators satisfy the same recursions. We consider two
types of isomorphisms in this paper: $\hilb=\ell^2(G/P)$, where $\Phi$ is
derived from the $\phi$ defined in~(\ref{eq:phi}); the second case
considered is $\hilb=L^2(\partial\tree)$, where
$\Phi:\hilb\to\hilb^\Sigma$ is defined by
$\Phi(f)(\sigma)=(f(0\sigma),\dots,f((d-1)\sigma))$, for $f\in
L^2(\partial\tree)$ and $\sigma\in\partial\tree$. There are actually
uncountably many non-equivalent isomorphisms, giving uncountably many
non-equivalent representations of the same group, as we indicated just
before Subsection~\ref{subs:approx}.

\section{Computations of Finite Spectra}
Here we compute explicitly the spectra of the representations $\pi_n$
described in Section~\ref{sec:urho}. For our five examples, the general
principle will be the same: obtain recurrence relations on the
matrices of the representation and compute eigenvalues by recurrence.
We end each subsection with a computation of the spectral measure
$\tau_\Delta$.

\subsection{The group $G$}
Recall the finite quotient $G_n$ acts faithfully on $\{0,1\}^n$.  If
we denote by $a_n, b_n, c_n, d_n$ the permutation matrices of this
representation, we have
\begin{gather*}
  a_0 = b_0 = c_0 = d_0 = (1),\\
  a_n = \begin{pmatrix}0&1\\1&0\end{pmatrix},\qquad b_n = \begin{pmatrix}a_{n-1}&0\\0&c_{n-1}\end{pmatrix},\\
  c_n = \begin{pmatrix}a_{n-1}&0\\0&d_{n-1}\end{pmatrix},\qquad d_n = \begin{pmatrix}1&0\\0&b_{n-1}\end{pmatrix}.
\end{gather*}

The Hecke type operator of $\pi_n$ is
\[\Delta_n = a_n + b_n + c_n + d_n = \begin{pmatrix}2a_{n-1}+1&1\\1&\Delta_{n-1}-a_{n-1}\end{pmatrix},\]
and we wish to compute its spectrum.  We start by proving a slightly
stronger result: define
\[Q_n(\lambda,\mu) = \Delta_n - (\lambda+1)a_n - (\mu+1)\]
and
\begin{align*}
  \Phi_0 &= 2-\mu-\lambda,\\
  \Phi_1 &= 2-\mu+\lambda,\\
  \Phi_2 &= \mu^2-4-\lambda^2,\\
  \Phi_n &= \Phi_{n-1}^2 - 2(2\lambda)^{2^{n-2}}\qquad(n\ge3).
\end{align*}

\begin{lem}\label{lem:QG}
  For $n\ge2$ we have
  \[|Q_n(\lambda,\mu)| = (4-\mu^2)^{2^{n-2}}\left|Q_{n-1}\left(\frac{2\lambda^2}{4-\mu^2},\mu+\frac{\mu\lambda^2}{4-\mu^2}\right)\right|\qquad(n\ge2).\]
\end{lem}
\begin{proof}
  \begin{align*}
    |Q_n(\lambda,\mu)| &= \begin{vmatrix}2a_{n-1}-\mu&-\lambda\\-\lambda&\Delta_{n-1}-a_{n-1}-(1+\mu)\end{vmatrix}\\
    &= \begin{vmatrix}2a_{n-1}-\mu&-\lambda\\(2a_{n-1}-\mu)\frac{-\lambda(2a_{n-1}+\mu)}{4-\mu^2}&\Delta_{n-1}-a_{n-1}-(1+\mu)\end{vmatrix}\\
    &= \begin{vmatrix}2a_{n-1}-\mu&-\lambda\\0&\Delta_{n-1}-a_{n-1}-(1+\mu)-\frac{\lambda^2}{4-\mu^2}(2a_{n-1}+\mu)\end{vmatrix}\\
    &= |2a_{n-1}-\mu|\cdot\left|\Delta_{n-1}-\left(1+\frac{2\lambda^2}{4-\mu^2}\right)a_{n-1}-\left(1+\mu+\frac{\mu\lambda^2}{4-\mu^2}\right)\right|;
  \end{align*}
  now \[|2a_{n-1}-\mu| = \begin{vmatrix}-\mu&2\\2&-\mu\end{vmatrix}_{2^{n-1}} = (\mu^2-4)^{2^{n-2}}\]
  so the lemma follows.
\end{proof}

\begin{lem}\label{lemma:PhiG}
  For all $n$ we have
  \[|Q_n| = \Phi_0\Phi_1\cdots\Phi_n.\]
\end{lem}
\begin{proof}
  By direct computation,
  \[Q_0(\lambda,\mu) = (2-\mu-\lambda),\qquad Q_1(\lambda,\mu) = \begin{pmatrix}2-\mu&-\lambda\\-\lambda&2-\mu\end{pmatrix},\]
  so $|Q_0|=\Phi_0$ and $|Q_1|=\Phi_0\Phi_1$.

  Let us temporarily agree to write
  \begin{xalignat*}{2}
    \lambda'&= \frac{2\lambda^2}{4-\mu^2} & \Phi_i &= \Phi_i(\lambda,\mu)\\
    \mu'&=\mu+\frac{\mu\lambda^2}{4-\mu^2} & \Phi_i' &= \Phi_i(\lambda',\mu').
  \end{xalignat*}
  We will prove by recurrence that
  \begin{align*}
    (2-\mu)\Phi_0' &= \Phi_0\Phi_1,\\
    (2+\mu)\Phi_1' &= -\Phi_2,\\
    (4-\mu^2)\Phi_2' &= -\Phi_3,\\
    (4-\mu^2)^{2^{n-2}}\Phi_n' &= \Phi_{n+1}\qquad(n\ge3).
  \end{align*}
  Indeed
  \begin{align*}
    (2-\mu)\Phi_0' &= (2-\mu)(2-\mu) - \frac{\mu\lambda^2}{2+\mu} - \frac{2\lambda^2}{2+\mu}\\
    &= (2-\mu)^2-\lambda^2 = \Phi_0\Phi_1;\\
    (2+\mu)\Phi_1' &= (2+\mu)(2-\mu) - \frac{\mu\lambda^2}{2-\mu} + \frac{2\lambda^2}{2-\mu}\\
    &= 4-\mu^2+\lambda^2 = -\Phi_2;\\
    (4-\mu^2)\Phi_2' &= (4-\mu^2)(\mu^2-4) + 2\mu^2\lambda^2 + \frac{(\mu\lambda^2)^2}{4-\mu^2} - \frac{(2\lambda^2)^2}{4-\mu^2}\\
    &= -\left[(\mu^2-4)^2 - 2\mu^2\lambda^2 + \lambda^4\right]\\
    &= -\left[(\mu^2-4-\lambda^2)^2 - 8\lambda^2\right] = -\Phi_3;\\
\intertext{and for $n\ge3$,}
    (4-\mu^2)^{2^{n-2}}\Phi_n' &= (4-\mu^2)^{2^{n-2}}(\Phi_{n-1}')^2-(4-\mu^2)^{2^{n-2}}2(2\lambda')^{2^{n-2}}\\
    &= \left((4-\mu^2)^{2^{n-3}}\Phi_{n-1}'\right)^2-2(4\lambda^2)^{2^{n-2}}\\
    &= (\pm\Phi_n)^2-2(2\lambda)^{2^{n-1}} = \Phi_{n+1}.
  \end{align*}
  Now using Lemma~\ref{lem:QG} we have for $n\ge3$
  \begin{align*}
    |Q_n(\lambda,\mu)| &= (4-\mu^2)^{2^{n-2}}|Q_{n-1}(\lambda',\mu')|\\
    &= (4-\mu^2)^{2^{n-2}}\Phi_0'\Phi_1'\cdots\Phi_{n-1}'\\
    &= (2-\mu)\Phi_0'(2+\mu)\Phi_1'(4-\mu^2)\Phi_2'\cdots(4-\mu)^{2^{n-3}}\Phi_{n-1}'\\
    &= \Phi_0\Phi_1\Phi_2\cdots\Phi_n,
  \end{align*}
  proving the claim.
\end{proof}  

\begin{prop}\label{prop:G}
  For all $n$ we have
  \begin{multline*}
    \{(\lambda,\mu):\,Q_n(\lambda,\mu)\text{ non invertible}\} = \{(\lambda,\mu):\,\Phi_0(\lambda,\mu)=0\}\cup \{(\lambda,\mu):\,\Phi_1(\lambda,\mu)=0\}\\
    {}\cup\{(\lambda,\mu):\,4-\mu^2+\lambda^2+4\lambda\cos\left(\frac{2\pi j}{2^n}\right)=0\text{ for some }j = 1,\dots,2^{n-1}-1\}.
  \end{multline*}
\end{prop}
\begin{proof}
  We prove by recurrence that for all $n,k$ with $0\le k\le n-2$ we
  have
  \begin{equation}\label{prop:pf1}
    \Phi_n = \prod_{t=0}^{2^k-1}\left(\Phi_{n-k} - 2(2\lambda)^{2^{n-2-k}}\cos\left(\frac{2\pi(2t+1)}{2^{k+2}}\right)\right).
  \end{equation}
  Indeed for $k=0$ equality holds trivially, and if $k>0$ we combine
  the terms for $t$ and $2^k-1-t$, with $t<2^{k-1}$: letting $A_{k,t}$ designate the $t$-th term,
  \begin{align*}
    A_{k,t}A_{k,2^k-1-t} &= \Phi_{n-k}^2 - 4(2\lambda)^{2^{n-k-1}}\cos^2\left(\frac{2\pi(2t+1)}{2^{k+2}}\right)\\
      &= \Phi_{n-k}^2 - 2\left(1+(2\lambda)^{2^{n-k-1}}\cos\left(\frac{2\pi(2t+1)}{2^{k+1}}\right)\right)\\
      &= \Phi_{n-k+1} - 2(2\lambda)^{2^{n-k-1}}\cos\left(\frac{2\pi(2t+1)}{2^{k+1}}\right)\\
      &= A_{k-1,t},
  \end{align*}
  so $\prod_{t=0}^{2^k-1} A_{k,t} = \prod_{t=0}^{2^{k-1}-1}A_{k-1,t}$.
  Letting $k=n-2$ in~(\ref{prop:pf1}) proves the proposition, in light
  of Lemma~\ref{lemma:PhiG}.
\end{proof}

In the $(\lambda,\mu)$ system, the spectrum of $Q_n$ is thus a
collection of $2$ lines and $2^{n-1}-1$ hyperbol\ae.
\begin{figure}
\begin{center}
\setlength\unitlength{1pt}
\begin{picture}(360,300)
\put(180,270){\makebox(0,0)[lb]{\smash{\normalsize $\mu$}}}
\put(310,130){\makebox(0,0)[lb]{\smash{\normalsize $\lambda$}}}
\put(0,-40){\epsfig{file=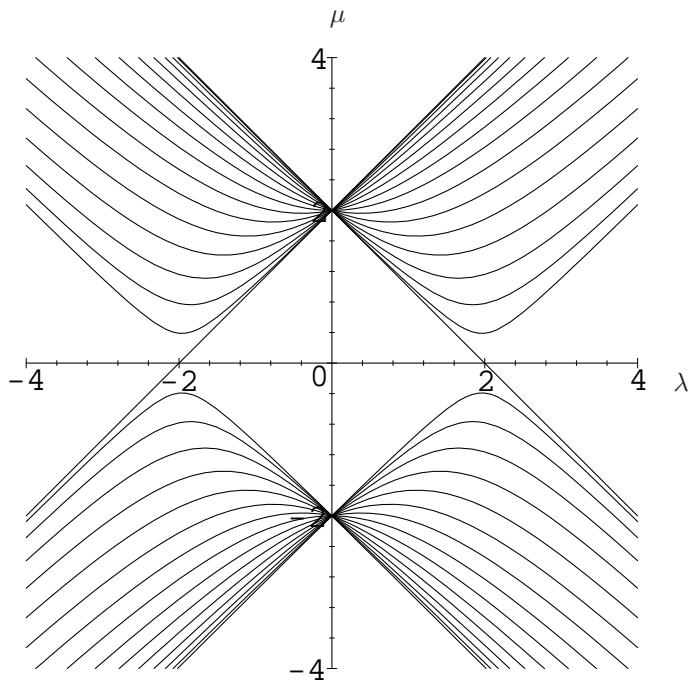}}
\end{picture}
\end{center}
\caption{The spectrum of $Q_n(\lambda,\mu)$ for $G$ and $\tilde G$}
\end{figure}
The spectrum of $\Delta_n$ is precisely the set of $\theta$ such that
$|Q_n(-1,\theta-1)|=0$. From Proposition~\ref{prop:G} we obtain
\[\spec(\Delta_n) = \{1\pm\sqrt{5-4\cos\phi}:\,\phi\in2\pi\Z/2^n\}\setminus\{0,-2\}.\]
The first eigenvalues are $4$; $2$; $1\pm\sqrt5$;
$1\pm\sqrt{5\pm2\sqrt2}$; $1\pm\sqrt{5\pm2\sqrt{2\pm\sqrt2}}$; etc.

The numbers of the form
$\pm\sqrt{\lambda\pm\sqrt{\lambda\pm\sqrt{\dots}}}$ appear as preimages
of the quadratic map $z^2-\lambda$, and after closure produce a Julia
set for this map (see~\cite{barnsley:fe}).  In the given example this
Julia set is just an interval. In the examples that follow in
Subsection~\ref{subs:computeGamma} the spectra are simple
transformations of Julia sets which are totally disconnected---this
behaviour is explained by Lemma~\ref{lem:dyn} and the remark after
its proof.
\begin{cor}
  The spectrum of $\pi$, for the group $G$, is
  \[\spec(\Delta) = [-2,0]\cup[2,4].\]
\end{cor}

We now investigate the empiric spectral measure $\tau_\Delta$, as
defined in Subsection~\ref{subs:sm}. We constructed in the previous
paragraph a one-to-one map $\chi:[0,\pi]\times\{\pm1\}\to\R$ defined
by
\[\chi(\theta,\epsilon)=1+\epsilon\sqrt{5+4\cos\theta}.\]
The spectrum is uniformly distributed in $[0,\pi]\times\{\pm1\}$ by
Proposition~\ref{prop:G}, and
$\chi$ is by assumption a measure-preserving map, so the measure
$\tau(A)$ of any subset $A$ of $\R$ can be evaluated as
\[\tau(A) = \operatorname{vol}(\chi^{-1}(A)),\]
where $\operatorname{vol}$ is the uniform measure on
$[0,\pi]\times\{\pm1\}$. The measure $d\tau_\Delta(x)=g(x)dx$ we are
seeking is thus given by
\[g(x) = \frac1{2\pi}\frac{\partial}{\partial x}\chi^{-1}(x) = \frac{1-x}{4\pi\sqrt{1-\left(\frac{(x-1)^2-5}{4}\right)^2}}.\]
\begin{figure}
\begin{center}
\setlength\unitlength{1pt}
\begin{picture}(360,300)
\put(180,270){\makebox(0,0)[lb]{\smash{\normalsize $f$}}}
\put(310,130){\makebox(0,0)[lb]{\smash{\normalsize $\lambda$}}}
\put(0,-40){\epsfig{file=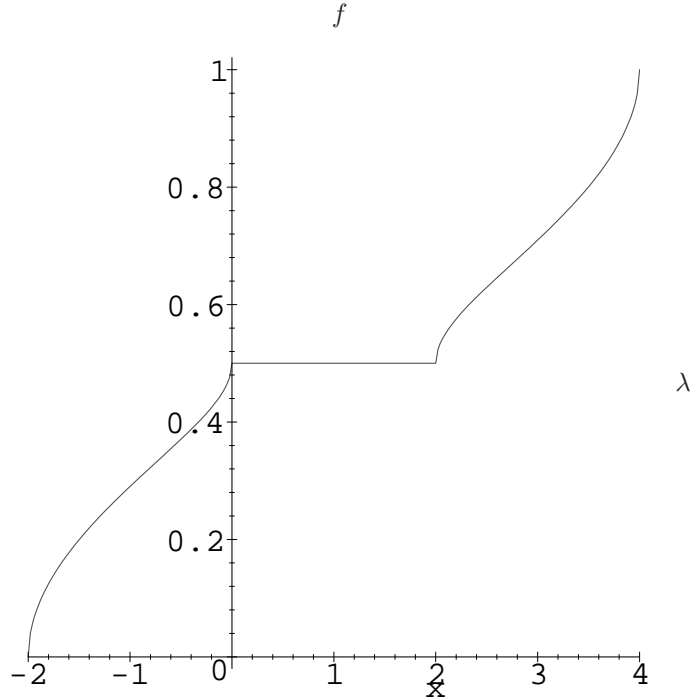}}
\end{picture}
\end{center}
\caption{The empiric spectral measure for $G$}
\end{figure}

The eigenvectors of $\Delta_n$ can be expressed as follows. Let
$\lambda$ be an eigenvalue and define by induction the sequence
\[v_1 = 1;\quad v_2=\lambda-3;\quad v_i=\begin{cases}\frac{(\lambda-1)v_{i-1}-v_{i-2}}{2}&\text{ if }i\ge3,i\equiv1[2];\\
  (\lambda-1)v_{i-1}-2v_{i-2}&\text{ if }i\ge3,i\equiv0[2].\end{cases}\]

Define also by induction for all $i\ge0$ the following ordering of
$\Sigma^i$: if the ordering of $\Sigma^{i-1}$ is
$(\sigma_1,\dots,\sigma_{2^{i-1}})$, the ordering of $\Sigma^i$ is
\[(1\sigma_1,0\sigma_1,0\sigma_2,1\sigma_2,1\sigma_3,\dots,0\sigma_{2^{i-1}},1\sigma_{2^{i-1}}).\]
\begin{lem}
  If $\lambda$ is an eigenvalue of $\Delta_n$, and
  $(\sigma_1,\dots,\sigma_{2^n})$ is the ordering described above,
  then the eigenvector corresponding to $\lambda$ has value $v_i$ on
  the vertex $\sigma_i$.
\end{lem}
\begin{proof}
  Consider the Schreier graph $\gf_n=\sch(G,P_n,S)$ described in
  Subsection~\ref{subs:schG}. It has vertex set $\{1,\dots,2^n\}$, a
  simple edge between $2i-1$ and $2i$ for all $i$, a double edge
  between $2i$ and $2i+1$ for all $i$, a loop at each $i$, and a
  triple loop at $1$ and $2^n$. The homogeneous space $G_n/P_n$ is
  isomorphic to $X$ graph through the correspondence $\sigma_i\mapsto
  i$. The choices of $v_i$ clearly define an eigenvector for
  $\lambda$, because $v_{2^n}=0$ if $\lambda$ is an eigenvalue for
  $\Delta_n$.
\end{proof}

It is then possible to express the characteristic function of any
vertex as a linear combination of the eigenvectors described above,
and to explicitly compute the corresponding Kesten spectral measure.
We shall not develop here this topic; but we note that the coordinates
of the eigenvectors (of norm $1$) vary greatly in absolute value, so
the spectral measure may be very different to the empiric spectral
measure.

\subsection{The Group $\tilde G$}\label{subs:spgt}
The computation of the spectrum of $\pi_n$ for $\tilde G$ amounts to
the following proposition. As before, we define $\tilde\Delta_n =
a_n+\tilde b_n+\tilde c_n+\tilde d_n$, viewed as a $2^n\times2^n$
matrix, and let
\[\tilde Q_n(\lambda,\mu) = \tilde\Delta_n - (\lambda+1)a_n-(\mu+2).\]
\begin{prop}
  Then for all $n$ we have
  \[\tilde Q_n(\lambda,\mu) = \frac12 Q_n(2\lambda,2\mu).\]
\end{prop}
\begin{proof}
  The proof follows by recurrence on $n$: first it is readily
  checked that $\tilde Q_0(\lambda,\mu) = (1-\mu-\lambda)$; then for
  $n\ge1$ we compute
  \begin{align*}
    \frac12 Q_n(2\lambda,2\mu) &= \frac12\begin{pmatrix}2a_{n-1}-2\mu&-2\lambda\\-2\lambda&\Delta_{n-1}-a_{n-1}-(2\mu+1)\end{pmatrix}\\
    &= \begin{pmatrix}a_{n-1}-\mu&-\lambda\\-\lambda&\frac12 Q_{n-1}(0,2\mu)\end{pmatrix}=\begin{pmatrix}a_{n-1}-\mu&-\lambda\\-\lambda&\tilde Q_{n-1}(0,\mu)\end{pmatrix}\\
    &= \begin{pmatrix}a_{n-1}-\mu&-\lambda\\-\lambda&\tilde\Delta_{n-1}-a_{n-1}-(\mu+2)\end{pmatrix}=\tilde Q_n(\lambda,\mu).
  \end{align*}
\end{proof}

We can now obtain the spectrum of $\pi_n$ by setting $\lambda=-1$
in $\tilde Q_n(\lambda,\mu)$ and solving for $\mu$; in view of the
previous proposition and the computations performed for $G_n$, we obtain:
\begin{prop}
  The spectrum of $\tilde G_n$ is
  \[\left.\left\{2+2\cos\left(\frac{2\pi j}{2^{n+1}}\right)\right|j=0,\dots,2^n-1\right\}.\]
\end{prop}
\begin{proof}
  The spectrum consists precisely of the $\mu+2$ such that $\tilde
  Q_n(-1,\mu)=0$; this amounts to $Q_n(-2,2\mu)=0$. Now this holds
  when $\Phi_0(-2,2\mu)=0$, $\Phi_1(-2,2\mu)=0$ or
  $4-(2\mu)^2+4-4\cos(2\pi j/2^n)(-2)=0$ for some
  $j=1,\dots,2^{n-1}-1$. These give respectively $\mu=2$, $\mu=0$ and
  $\mu=\pm\sqrt{2-2\cos(2\pi j/2^n)}$, which after simplification
  yield the proposition.
\end{proof}
The first eigenvalues are $4$; $2$; $2\pm\sqrt2$;
$2\pm\sqrt{2\pm\sqrt2}$; $2\pm\sqrt{2\pm\sqrt{2\pm\sqrt2}}$; etc.

\begin{cor}
  The spectrum of $\pi$, for the group $\tilde G$, is
  \[\spec(\Delta) = [0,4].\]
\end{cor}
Note that the spectrum of $\Delta$ is positive! This can never happen
to regular representations~\cite{harpe-r-v:sg}, and indeed
the spectrum of the regular representation of $\tilde G$ is
$[-4,4]$, since it contains $[0,4]$ and is symmetrical about $0$ (as
the Cayley graph $\mathcal C(\tilde G,S)$ is bipartite).

Again we may compute the empiric spectral measure $\tau_\Delta$.
Recall that there is a measure-preserving map $\chi:[0,\pi]\to\R$
given by $\theta\mapsto2+2\cos\theta$, where $[0,\pi]$ has the
Lebesgue measure. The measure $d\tau_\Delta=g(x)dx$ we are seeking is then
associated to
\[g(x) = \frac1\pi\frac{\partial}{\partial x}\chi^{-1}(x) = \frac1{\pi\sqrt{4x-x^2}}.\]

\subsection{The Group $\Gamma$}\label{subs:computeGamma}
Recall that the finite quotient $\Gamma_n$ acts faithfully on
$\{0,1,2\}^n$. Denote by $a_n$ and $s_n$ the matrices of the action.
We have
\begin{gather*}
  a_0 = s_0 = (1),\\
  a_n = \begin{pmatrix}0&1&0\\0&0&1\\1&0&0\end{pmatrix},\qquad s_n = \begin{pmatrix}a_{n-1}&0&0\\0&1&0\\0&0&s_{n-1}\end{pmatrix}.
\end{gather*}

The combinatorial Laplacian of $\Gamma_n$ is
\[\Delta_n = a_n + a_n^{-1} + s_n + s_n^{-1} =
  \begin{pmatrix}a_{n-1}+a_{n-1}^{-1}&1&1\\1&a_{n-1}+a_{n-1}^{-1}&1\\1&1&s_{n-1}+s_{n-1}^{-1}\end{pmatrix}.\]
For ease of notation, let us define operators
\begin{gather*}
  A_n = a_n + a_n^{-1},\qquad S_n = s_n + s_n^{-1},\\
  Q_n(\lambda,\mu) = S_n + \lambda A_n - \mu,\\
  \intertext{and polynomials}
  \alpha=2-\mu+\lambda,\qquad \beta=2-\mu-\lambda,\\
  \gamma=\mu^2-\lambda^2-\mu-2,\qquad\delta=\mu^2-\lambda^2-2\mu-\lambda.
\end{gather*}

\begin{lem}\label{lem:GF1}
  We have
  \begin{gather*}
    |Q_0(\lambda,\mu)| = 2+2\lambda-\mu=\alpha+\lambda,\\
    |Q_1(\lambda,\mu)| = (2+2\lambda-\mu)(2-\lambda-\mu)^2=(\alpha+\lambda)\beta^2,\\
    |Q_n(\lambda,\mu)| = (\alpha\beta\gamma^2)^{3^{n-2}}\left|Q_{n-1}\left(\frac{\lambda^2\beta}{\alpha\gamma},\mu+\frac{2\lambda^2\delta}{\alpha\gamma}\right)\right|\quad(n\ge2).
  \end{gather*}
\end{lem}
\begin{proof}
  We compute the determinant:
  \begin{align*}
    |Q_n(\lambda,\mu)| &= |S_n + \lambda A_n - \mu| =
    \begin{vmatrix}A_{n-1}-\mu&\lambda&\lambda\\
      \lambda&2-\mu&\lambda\\ \lambda&\lambda&S_{n-1}-\mu\end{vmatrix}\\
    &= \begin{vmatrix}2-\mu-\lambda&\lambda+\mu-A_{n-1}&0\\
      \lambda&A_{n-1}-\mu&\lambda\\
      \lambda&\lambda&S_{n-1}-\mu\end{vmatrix}\\
    &= \begin{vmatrix}1&0&0\\
      *&(2-\mu)A_{n-1}+(\mu^2-\lambda^2-2\mu)&\lambda\\
      *&\lambda A_{n-1}+\lambda(2-2\lambda-2\mu)&S_{n-1}-\mu\end{vmatrix}\\
    &= \left|S_{n-1}-\mu-\frac{\lambda^2(1-2\lambda-\mu+A_{n-1})}{(1-\mu)A_{n-1}+\mu^2-\lambda^2}\right|\cdot|(1-\mu)A_{n-1}+\mu^2-\lambda^2|\\
    &= \left|S_{n-1}-\mu+\lambda^2\frac{\beta A_{n-1}+2\delta}{\alpha\gamma}\right|(\alpha\beta\gamma^2)^{3^{n-2}},
  \end{align*}
  which completes the proof of the lemma, using the easily verified
  equations
  \begin{gather*}
    \frac1{\pi A_n+\rho} = \frac{\pi A_n-\pi-\rho}{(2\pi+\rho)(\pi-\rho)},\\
    |\pi A_n+\rho| = \begin{vmatrix}\rho&\pi&\pi\\\pi&\rho&\pi\\\pi&\pi&\rho\end{vmatrix}_{3^n} = \left((\pi-\rho)^2(2\pi+\rho)\right)^{3^{n-1}}
  \end{gather*}
  valid for all scalar $\pi$ and $\rho$.
\end{proof}

\begin{figure}
  \begin{picture}(400,200)(-80,90)
    \put(0,0){\epsfig{file=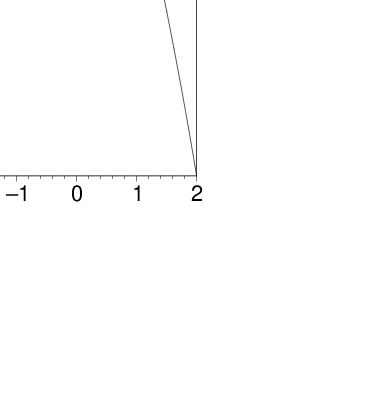}}
    \put(200,0){\epsfig{file=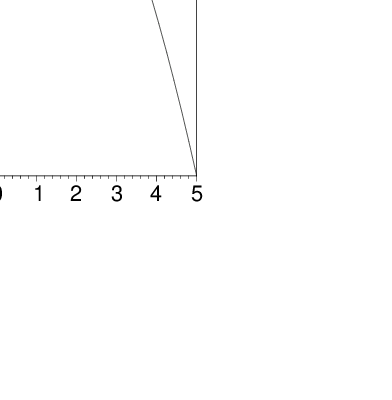}}
  \end{picture}
  \caption{The Functions $F$ for $\Gamma$ and $\overline\Gamma$}
  \label{fig:F}
\end{figure}
Consider now the quadratic forms
\[H_\theta = \mu^2-\lambda\mu-2\lambda^2-2-\mu+\theta\lambda,\]
and the function $F:[-4,5]\to[-4,5]$ (see Figure~\ref{fig:F} left)
given by
\[F(\theta) = 4-2\theta-\theta^2.\]
Let $X_2 = \{-1\}$, and iteratively define $X_n = F^{-1}(X_{n-1})$
for all $n\ge3$.  Note that $|X_n| = 2^{n-2}$.
\begin{lem}\label{lem:factorQGF}
  We have for all $n\ge 2$ the factorization
  \begin{equation}
    |Q_n(\lambda,\mu)|=(2+2\lambda-\mu)(2-\lambda-\mu)^{3^{n-1}+1}\prod_{\substack{2\le m\le n\\ \theta\in X_m}}H_\theta^{3^{n-m}+1}.
  \end{equation}
\end{lem}
\begin{proof}
  Let us write temporarily
  \[\lambda' = \frac{\lambda^2\beta}{\alpha\gamma},\qquad\mu' = \mu+\frac{2\lambda^2\delta}{\alpha\gamma},\]
  and $P'=P(\lambda',\mu')$ for $P$ in
  $\{\alpha,\beta,\gamma,\delta,H_\theta\}$. Then we have
  \begin{gather*}
    \alpha'+\lambda'=\frac{\beta(\alpha+\lambda)}{\alpha},\qquad\beta'=\frac{\beta H_{-1}}{\gamma},\\
    H_\theta(\lambda',\mu') = \frac{\beta\prod_{\nu\in F^{-1}(\theta)}H_\nu(\lambda,\mu)}{\alpha\gamma}.
  \end{gather*}
  Now, first $Q_2(\lambda,\mu)=(\alpha+\lambda)\beta^4H_{-1}^2$ as
  claimed; then by induction, using Lemma~\ref{lem:GF1}, we have for $n\ge3$
  \begin{align*}
    \begin{split}
      |Q_n(\lambda,\mu)| &= (\alpha\beta\gamma2)^{3^{n-2}}\left|Q_{n-1}(\lambda',\mu')\right|\\
      &= (\alpha\beta\gamma^2)^{3^{n-2}}(\alpha'+\lambda'){\beta'}^{3^{n-2}+1}\prod_{\substack{2\le m\le n-1\\ \theta\in X_m}}{H'_\theta}^{3^{n-1-m}+1}\\
      &= (\alpha\beta\gamma^2)^{3^{n-2}}\frac{\beta(\alpha+\lambda)}{\alpha}\left(\frac{\beta H_{-1}}{\gamma}\right)^{3^{n-2}+1}\left(\frac{\beta}{\alpha\gamma}\right)^{3^{n-2}-1}\prod_{\substack{2\le m\le n\\ \theta\in X_m}}{H_\theta}^{3^{n-m}+1}\\
      &= (\alpha+\lambda)\beta^{3^{n-1}+1}\prod_{\substack{2\le m\le n\\ \theta\in X_m}}{H_\theta}^{3^{n-m}+1}
    \end{split}
  \end{align*}
  as claimed.
\end{proof}

Thus, according to the previous proposition, the spectrum of $Q_n$ is
a collection of two lines and $2^{n-1}-1$ hyperbol\ae\ $H_\theta$ with
$\theta\in X_2\sqcup X_3\sqcup\dots\sqcup X_n$. The spectrum of
$\Delta_n$ is obtained by solving $|Q_n(1,\mu)|=0$, as given in the
following proposition.
\begin{figure}
\begin{center}
\setlength\unitlength{1pt}
\begin{picture}(360,300)
\put(150,270){\makebox(0,0)[lb]{\smash{\normalsize $\mu$}}}
\put(310,100){\makebox(0,0)[lb]{\smash{\normalsize $\lambda$}}}
\put(0,-40){\epsfig{file=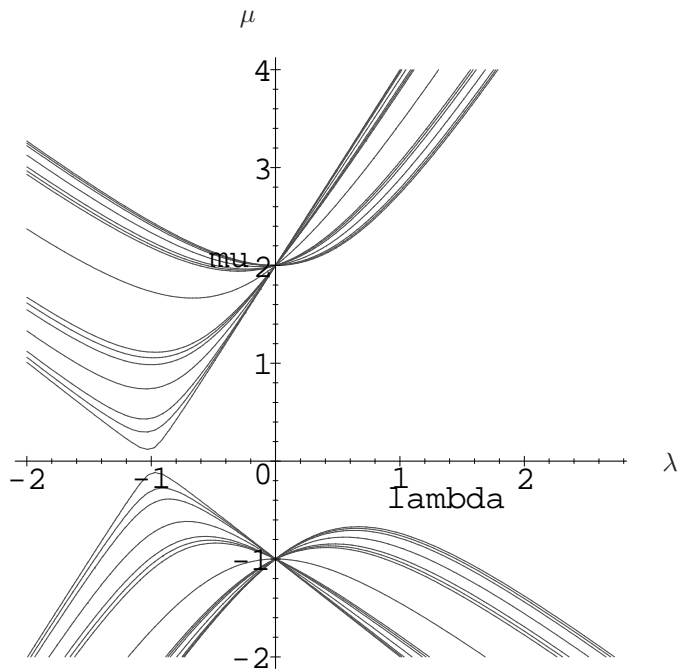}}
\end{picture}
\end{center}
\caption{The spectrum of $Q_n(\lambda,\mu)$ for $\Gamma$}
\end{figure}

\begin{prop}\label{prop:spGF}
  Let $\pi_\pm:[-4,5]\to[-2,4]$ be defined by
  $\pi_\pm(\theta)=1\pm\sqrt{5-\theta}$. Then
  \begin{gather*}
    \spec{\Delta_0} = \{4\};\\
    \spec{\Delta_1} = \{1,4\};\\
    \spec{\Delta_n} = \{1,4\} \cup \bigcup_{2\le m\le n}\pi_\pm(X_m)\qquad(n\ge 2).
  \end{gather*}
\end{prop}
\begin{proof}
  Solving $H_\theta(1,\mu)=0$ gives $\mu=\pi_\pm(\theta)$. The result
  then follows from the factorization of $|Q_n|$ as a product of
  $H_\theta$'s given by Lemma~\ref{lem:factorQBG}.
\end{proof}

\begin{figure}
\epsfig{file=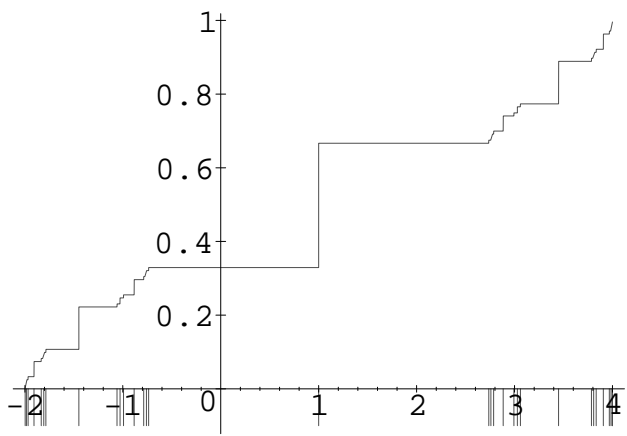}
\caption{The Cantor-Set spectrum of $\Delta_n$ for $\Gamma$}
\end{figure}

\begin{lem}\label{lem:dyn}
  Let $F:[a,b]\to\R$ be a smooth unimodal map with negative Schwartzian
  derivative (for instance, a quadratic map), and $F(a)=F(b)=a$. Choose
  $\xi\in[a,b]$, and form
  \[K=\bigcap_{n=0}^\infty F^{-n}([a,b]),\qquad
  L=\overline{\bigcup_{n=0}^\infty F^{-n}(\{\xi\})}.\]
  Then the following cases may occur:
  \begin{description}
  \item[$\mathbf{F([a,b])\subset[a,b]}$] Then $K=[a,b]$.
  \item[$\mathbf{F([a,b])=[a,b]}$] Then $K=L=[a,b]$.
  \item[$\mathbf{F([a,b])\supset[a,b]}$] Then $K$ is a Cantor set of
    null Lebesgue measure. If $\xi\in K$, then $L=K$, otherwise
    $L=K\sqcup C$, where $C$ is a countable set of isolated points
    accumulating on $K$.
  \end{description}
\end{lem}
\begin{proof}
  In the first case: $K$ is clearly $[a,b]$, because
  $F^{-1}([a,b])=[a,b]$.
  
  Now assume that $F([a,b])\supseteq[a,b]$, and define two continuous
  maps $g_1,g_2:[a,b]\to[a,b]$ such that
  $F^{-1}(x)=\{g_1(x),g_2(x)\}$.  Then $\Pi=\langle g_1,g_2\rangle$ is
  a free semigroup, because $g_1$ and $g_2$ have disjoint
  image-interiors.
  
  In the second case: $K$ is obviously all $[a,b]$, because
  $F^{-1}([a,b])=[a,b]$. Clearly $\Pi$-orbit of $\xi$, is dense in
  $[a,b]=g_1([a,b])\cup g_2([a,b])$, so its closure $L$ is $[a,b]$.
  
  In the last case: we first show that $K\subseteq L$. Since
  $K=g_1(K)\sqcup g_2(K)$, every point in $K$ is specified by an
  infinite sequence of $g_i$'s applied to some point. Since the
  $g_i$'s are contracting, the choice of that point is unimportant ---
  hence we may choose $\xi$, and this proves the claim. If $\xi\in K$,
  then clearly $L\subseteq K$ and we are done. Otherwise, set $C=\cup
  F^{-n}(\{\xi\})$, which is a countable set of isolated points.
  Clearly $K$ and $C$ are disjoint (else $\xi$ would be in $K$).
  Finally $C$ accumulates on the infinite orbits under $\langle
  g_1,g_2\rangle$, hence on $K$.
  
  Let $O=[a,b]\setminus F^{-1}([a,b])$. Then $O$ disconnects $[a,b]$,
  and $F^{-1}(O)$ disconnects each of the connected components of
  $[a,b]\setminus O$, etc., so $K$ is disconnected. Take any point
  $x\in K$; then for all $n\in\N$ we may write $x=w(x_n)$ for some
  $x_n\in K$ and some word $w\in\Pi$ of length $n$. Taking some
  $x'_n\in K$ close to $x_n$ and letting $n$ tend to infinity gives a
  sequence in $K$ converging to $x$. Since $K$ is clearly closed, it
  is perfect, so is a Cantor set.  Finally
  $m(K)=\lim_{n\to\infty}(1-m(O))^n=0$, where $m$ denotes the Lebesgue
  measure.
\end{proof}
Let us make a few remarks concerning this lemma:
\begin{itemize}
\item The behaviour of $L$ in the first case seems unpredictable. It
  is fortunately not needed for our purposes.
\item The second case arises in connection with $G$ and $\tilde G$,
  the third (with $\xi\not\in K$) with $\Gamma$ and (with $\xi\in K$)
  with $\overline\Gamma$ and $\doverline\Gamma$.
\item Much more is known on the forward and inverse orbits of points
  under a unimodal map. See for instance~\cite[pages~10
  and~327--351]{demelo-s:1dd} for more information.
\item $K$ is the Julia set of the dynamical system $F$, i.e.\ the
  closure of the set of points whose (forward) orbit does not diverge
  to infinity. Indeed $K$ consists precisely in those points whose
  orbit remains in $[a,b]$.  Examples of Julia sets of a similar
  nature (occurring in the study of the quadratic map), with nested
  square-root expressions, appear in~\cite[page~277]{barnsley:fe} ---
  see also the bibliography there.
\end{itemize}

The previous lemma gives the structure of the spectrum of $\Gamma$:
\begin{cor}
  The spectrum of $\pi$ for the group $\Gamma$ is the closure of the set
  \[\left\{\begin{array}{rl}
      4 & (\mu=0)\\
      1 & (\mu=\frac13)\\
      1\pm\sqrt6 & (\mu=\frac19)\\
      1\pm\sqrt{6\pm\sqrt6} & (\mu=\frac1{3^3})\\
      1\pm\sqrt{6\pm\sqrt{6\pm\sqrt6}} & (\mu=\frac1{3^4})\\
      \dots\end{array}\right\}.\]
  It is the union of a Cantor set of
  null Lebesgue measure that is symmetrical about $1$, and a countable
  collection of isolated points supporting the empiric spectral
  measure, which has the values indicated as $\mu$.
\end{cor}

\subsection{The Group $\overline\Gamma$}
Recall the finite quotient $\Gamma_n$ acts faithfully on
$\{0,1,2\}^n$. Denote by $a_n$ and $t_n$ the matrices of this action.
We have
\begin{gather*}
  a_0 = t_0 = (1),\\
  a_n = \begin{pmatrix}0&0&1\\1&0&0\\0&1&0\end{pmatrix},\qquad
  t_n = \begin{pmatrix}a_{n-1}&0&0\\0&a_{n-1}&0\\0&0&t_{n-1}\end{pmatrix}.
\end{gather*}

The combinatorial Laplacian of $\Gamma_n$ is
\[\Delta_n = a_n + a_n^{-1} + t_n + t_n^{-1} =
  \begin{pmatrix}a_{n-1}+a_{n-1}^{-1}&1&1\\1&a_{n-1}+a_{n-1}^{-1}&1\\1&1&t_{n-1}+t_{n-1}^{-1}\end{pmatrix}.\]
For ease of notation, let us define operators
\begin{gather*}
  A_n = a_n + a_n^{-1},\qquad T_n = t_n + t_n^{-1},\\
  Q_n(\lambda,\mu) = T_n + \lambda A_n - \mu,\\
  \intertext{and polynomials}
  \alpha = 2-\mu+\lambda,\qquad\beta=2-\mu-\lambda,\\
  \gamma = 1+\mu+\lambda,\qquad\delta=1+\mu-\lambda.\\
\end{gather*}

\begin{lem}\label{lem:BG1}
  We have
  \begin{gather*}
    |Q_0(\lambda,\mu)| = \alpha+\lambda,\\
    |Q_1(\lambda,\mu)| = (\alpha+\lambda)\beta^2,\\
    |Q_n(\lambda,\mu)| = (\gamma\delta)^{2\cdot3^{n-2}}(\alpha\beta)^{3^{n-2}}\left|Q_{n-1}\left(\frac{-2\lambda^2}{\alpha\delta},\mu+\frac{2\lambda^2(\mu-\lambda-1)}{\alpha\delta}\right)\right|\quad(n\ge2).
  \end{gather*}
\end{lem}
\begin{proof}
  We compute the determinant:
  \begin{align*}
    |Q_n(\lambda,\mu)| &= |T_n + \lambda A_n - \mu| = \begin{vmatrix}A_{n-1}-\mu&\lambda&\lambda\\ \lambda&A_{n-1}-\mu&\lambda\\ \lambda&\lambda&T_{n-1}-\mu\end{vmatrix}\\
    &= \begin{vmatrix}A_{n-1}-\mu-\lambda&0&0\\ \lambda&A_{n-1}-\mu+\lambda&\lambda\\ \lambda&2\lambda&T_{n-1}-\mu\end{vmatrix}\\
    &= \left|T_{n-1}-\mu-\frac{2\lambda^2}{A_{n-1}-\mu+\lambda}\right|\cdot|A_{n-1}-\mu+\lambda|\cdot|A_{n-1}-\mu-\lambda|\\
    &= \left|T_{n-1}-\mu-2\lambda^2\frac{A_{n-1}+\mu-\lambda-1}{(2-\mu+\lambda)(1+\mu-\lambda)}\right|\cdot|A_{n-1}-\mu+\lambda|\cdot|A_{n-1}-\mu-\lambda|
  \end{align*}
  which completes the proof of the lemma, using the easily
  verified equations
  \begin{gather*}
    \frac1{A_n-\theta} = \frac{A_n+\theta-1}{(2-\theta)(1+\theta)},\\
    |A_n-\theta| = \begin{vmatrix}-\theta&1&1\\1&-\theta&1\\1&1&-\theta\end{vmatrix}_{3^{n-1}} = \left((\theta+1)^2(2-\theta)\right)^{3^{n-1}}
  \end{gather*}
  valid for all scalar $\theta$.
\end{proof}

Consider again the quadratic forms
\[H_\theta = \mu^2-\lambda\mu-2\lambda^2-2-\mu+\theta\lambda,\]
and consider the real function (see Figure~\ref{fig:F} right)
\[F(\theta) = \frac{12+\theta-\theta^2}{2}.\]
Let $X_3 = \{2\}$, $Y_3 = \{-1\}$ and iteratively define $X_n =
F^{-1}(X_{n-1})$ and $Y_n = F^{-1}(Y_{n-1})$ for all $n\ge4$.
Note that $|X_n| = |Y_n| = 2^{n-3}$.
\begin{lem}\label{lem:factorQBG}
  We have for all $n\ge 2$ the factorization
  \begin{equation}
    |Q_n(\lambda,\mu)| = (\alpha+\lambda)\beta^{3^{n-2}+1}\gamma^{3^{n-1}-1}(\delta-\lambda)^{3^{n-2}-1}\prod_{\substack{3\le m\le n\\ \theta\in X_m}}H_\theta^{3^{n-m}+1}\prod_{\substack{3\le m<n\\ \theta\in Y_m}}H_\theta^{3^{n-m}-1}\prod_{\theta\in X_{n+1}}H_\theta^2.
  \end{equation}
\end{lem}
\begin{proof}
  Let us write temporarily
  \[\lambda' = \frac{-2\lambda^2}{\alpha\delta},\qquad\mu' = \mu+\frac{2\lambda^2(\mu-\lambda-1)}{\alpha\delta},\]
  and $P'=P(\lambda',\mu')$ for $P$ in
  $\{\alpha,\beta,\gamma,\delta,H_\theta\}$. Then we have
  \begin{gather*}
    \alpha'+\lambda'=\frac{\beta(\alpha+\lambda)}{\alpha},\qquad\beta'=\frac{-H_2}{\delta},\qquad\gamma'=\frac{\gamma(\delta-\lambda)}{\delta},\qquad\delta'-\lambda'=\frac{-H_{-1}}{\alpha},\\
    H_\theta(\lambda',\mu') = \frac{\prod_{\nu\in F^{-1}(\theta)}H_\nu(\lambda,\mu)}{\alpha\delta}.
  \end{gather*}
  Now, first $Q_2(\lambda,\mu)=(\alpha+\lambda)\beta^2\gamma^2H_2^2$ as
  claimed; then by induction, using Lemma~\ref{lem:BG1}, we have for $n\ge3$
  \begin{align*}
    \begin{split}
      |Q_n(\lambda,\mu)| &= (\gamma\delta)^{2\cdot3^{n-2}}(\alpha\beta)^{3^{n-2}}\left|Q_{n-1}(\lambda',\mu')\right|\\
      &= (\gamma\delta)^{2\cdot3^{n-2}}(\alpha\beta)^{3^{n-2}}(\alpha'+\lambda'){\beta'}^{3^{n-3}+1}{\gamma'}^{3^{n-2}-1}(\delta'-\lambda')^{3^{n-3}-1}\\ &\hspace{6em}\prod_{\substack{3\le m\le n-1\\ \theta\in X_m}}{H'_\theta}^{3^{n-1-m}+1}\prod_{\substack{3\le m<n-1\\ \theta\in Y_m}}{H'_\theta}^{3^{n-1-m}-1}\prod_{\theta\in X_n}{H'_\theta}^2\\
      &= (\gamma\delta)^{2\cdot3^{n-2}}(\alpha\beta)^{3^{n-2}}\frac{\beta(\alpha+\lambda)}{\alpha}\left(\frac{-H_2}{\delta}\right)^{3^{n-3}+1}\left(\frac{\gamma(\delta-\lambda)}{\delta}\right)^{3^{n-2}-1}\left(\frac{-H_{-1}}{\alpha}\right)^{3^{n-3}-1}\\ &\hspace{6em}(\alpha\delta)^{-2\cdot 3^{n-3}}\prod_{\substack{4\le m\le n\\ \theta\in X_m}}H_\theta^{3^{n-m}+1}\prod_{\substack{4\le m<n\\ \theta\in Y_m}}H_\theta^{3^{n-m}-1}\prod_{\theta\in X_{n+1}}H_\theta^2.\\
      &= (\alpha+\lambda)\beta^{3^{n-2}+1}\gamma^{3^{n-1}-1}(\delta-\lambda)^{3^{n-2}-1}\prod_{\substack{3\le m\le n\\ \theta\in X_m}}H_\theta^{3^{n-m}+1}\prod_{\substack{3\le m<n\\ \theta\in Y_m}}H_\theta^{3^{n-m}-1}\prod_{\theta\in X_{n+1}}H_\theta^2
    \end{split}
  \end{align*}
  as claimed.
\end{proof}

Thus, according to the previous proposition, the spectrum of $Q_n$ is
a collection of $4$ lines and $(2^{n-1}-1)+(2^{n-3}-1)$ hyperbol\ae\
$H_\theta$ with $\theta\in X_3\sqcup\dots\sqcup X_{n+1}\sqcup Y_3\sqcup
\dots\sqcup Y_{n-1}$. The spectrum of $\Delta_n$ is obtained by
solving $|Q_n(1,\mu)|=0$, as given in the following proposition.
\begin{figure}
\begin{center}
\setlength\unitlength{1pt}
\begin{picture}(360,300)
\put(150,270){\makebox(0,0)[lb]{\smash{\normalsize $\mu$}}}
\put(310,100){\makebox(0,0)[lb]{\smash{\normalsize $\lambda$}}}
\put(0,-40){\epsfig{file=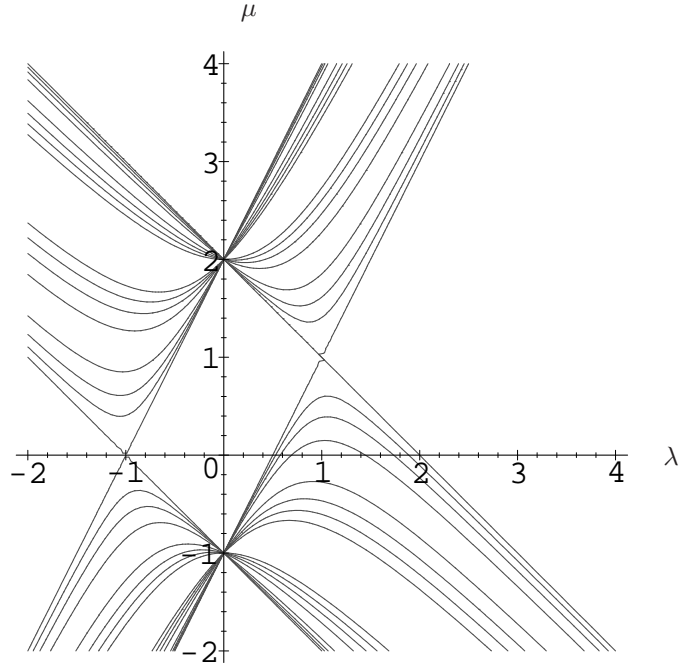}}
\end{picture}
\end{center}
\caption{The spectrum of $Q_n(\lambda,\mu)$ for $\overline\Gamma$ and
  $\doverline\Gamma$}
\end{figure}

\begin{prop}
  Let $\pi_\pm:[-4,5]\to[-2,4]$ be defined by $\pi_\pm(\theta) = 1 \pm
  \sqrt{5-\theta}$. Then
  \begin{gather*}
    \spec{\Delta_0} = \{4\};\\
    \spec{\Delta_1} = \{1,4\};\\
    \spec{\Delta_n} = \{-2,1,4\} \cup \bigcup_{3\le m\le n+1}\pi_\pm(X_m)
      \cup \bigcup_{3\le m\le n-1}\pi_\pm(Y_m)\qquad(n\ge 2).
  \end{gather*}
\end{prop}

\begin{figure}
\epsfig{file=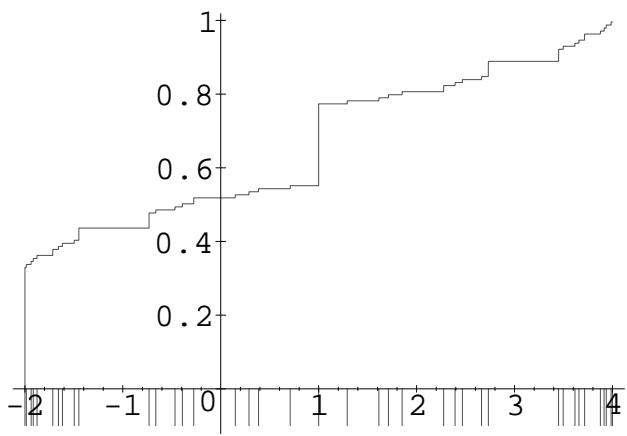}
\caption{The Cantor-Set spectrum of $\Delta_n$ for $\overline\Gamma$
  and $\doverline\Gamma$}
\end{figure}

\begin{proof}
  Solving $H_\theta(1,\mu)=0$ gives $\mu=\pi_\pm(\theta)$. The result
  then follows from the factorization of $|Q_n|$ as a product of
  $H_\theta$'s given by Lemma~\ref{lem:factorQBG}.
\end{proof}

\begin{cor}
  The spectrum of $\pi$ for the group $\overline\Gamma$ is the closure
  of the set
  \[\left\{\begin{array}{rl}4 & (\mu=0)\\
      -2 & (\mu=\frac13)\\
      1 & (\mu=\frac29)\\
      1\pm\sqrt{\frac{9\pm3}{2}} & (\mu=\frac2{27})\\
      1\pm\sqrt{\frac{9\pm\sqrt{45\pm4\cdot3}}{2}} & (\mu=\frac2{3^4})\\
      1\pm\sqrt{\frac{9\pm\sqrt{45\pm4\sqrt{45\pm4\cdot3}}}{2}} &
      (\mu=\frac2{3^5})\\
      \dots\end{array}\right\}.\]
  It is a Cantor set of null Lebesgue
  measure that is symmetrical about $1$. The empiric spectral measure
  is concentrated on the above algebraic numbers and has the values
  indicated as $\mu$.
\end{cor}

\subsection{The Group $\doverline\Gamma$}
Although $\overline\Gamma$ and $\doverline\Gamma$ greatly
differ in structure---the former is virtually torsion-free while the
latter is torsion---their representations $\pi_n$ have the same
spectrum:
\begin{prop}
  The spectrum of $\doverline\Gamma$ is the same as that of
  $\overline\Gamma$; that is, a Cantor set symmetrical about $1$ and
  spanning from $-2$ to $4$.
\end{prop}
\begin{proof}
  It suffices to note that for all $n$ we have
  $r_n+r_n^{-1}=t_n+t_n^{-1}$; this follows by induction on $n$.
\end{proof}

\section{Schreier Graphs}\label{sec:schreier}
Recall from Section~\ref{sec:urho} the definition of a Schreier graph.
In the general setting of a group $G$ acting level-transitively on a
tree $\Sigma^*$, and a subgroup $P$ defined as the stabilizer of an
infinite path, the quotient space $G/P$ is naturally identified with
an orbit, i.e.\ a countable subset, of $\Sigma^\N$. For the finite
quotients $G_n$, the space $G_n/P_n$ is identified with $\Sigma^n$. We
set $\gf_n$ denote the Schreier graph associated to this homogeneous
space. Due to the fractal (or recursive) nature of the five examples
we study, there are simple local rules producing $\gf_{n+1}$ from
$\gf_n$, the limit of this process being the Schreier graph of $G/P$.
We describe these rules for our examples: $G,\tilde
G,\Gamma,\overline\Gamma,\doverline\Gamma$.

\subsection{$\sch(G,P,S)$}\label{subs:schG}
We describe the graphs $\gf_n=\sch(G_n,P_n,S)$. They will be
with edges labeled by $S$ (and not oriented, because all $s\in S$ are
involutions) and vertices labeled by $\Sigma^n$, where
$\Sigma=\{0,1\}$.

First, it is clear that $\gf_0$ is a graph on one vertex, labeled by
the empty sequence $\emptyset$, and four loops at this vertex,
labeled by $a,b,c,d$. Next, $\gf_1$ has two vertices, labeled by $0$
and $1$; an edge labeled $a$ between them; and three loops at $0$ and
$1$ labeled by $b,c,d$.

Now given $\gf_n$, for some $n\ge1$, perform on it the following
transformation: replace the edge-labels $b$ by $d$, $d$ by $c$, $c$ by $b$;
replace the vertex-labels $\sigma$ by $1\sigma$; and replace all edges
labeled by $a$ connecting $\sigma$ and $\tau$ by: edges from
$1\sigma$ to $0\sigma$ and from $1\tau$ to $0\tau$, labeled $a$; two
edges from $0\sigma$ to $0\tau$ labeled $b$ and $c$; and loops at
$0\sigma$ and $0\tau$ labeled $d$. We claim the resulting graph is
$\gf_{n+1}$.

To prove the claim, it suffices to check that the letters on the
edge-labels act as described on the vertex-labels. If
$b(\sigma)=\tau$, then $d(1\sigma)=1\tau$, and similarly for $c$ and
$d$; this explains the cyclic permutation of the labels $b,c,d$. The
other substitutions are verified similarly.

As an illustration, here are $\gf_2$ and $\gf_3$ for $G$:
\begin{center}
  \begin{picture}(120,40)
    \put(0,20){\line(1,0){40}}\put(20,30){\msmash{a}}
    \put(40,20){\curve(0,0,20,3,40,0)\curve(0,0,20,-3,40,0)}
    \put(60,30){\msmash{b,c}}
    \put(80,20){\line(1,0){40}}\put(100,30){\msmash{a}}
    \put(0,15){\msmash{10}}
    \put(0,20){\spline(0,0)(18,27)(-18,27)(0,0)}\put(-25,35){\msmash{b,c,d}}
    \put(40,15){\msmash{00}}
    \put(40,20){\spline(0,0)(15,25)(-15,25)(0,0)}\put(40,35){\msmash{d}}
    \put(80,15){\msmash{01}}
    \put(80,20){\spline(0,0)(15,25)(-15,25)(0,0)}\put(80,35){\msmash{d}}
    \put(120,15){\msmash{11}}
    \put(120,20){\spline(0,0)(15,25)(-15,25)(0,0)}\put(142,35){\msmash{b,c,d}}
  \end{picture}\\
  \begin{picture}(280,40)
    \put(0,20){\line(1,0){40}}\put(20,30){\msmash{a}}
    \put(40,20){\curve(0,0,20,3,40,0)\curve(0,0,20,-3,40,0)}\put(60,30){\msmash{b,d}}
    \put(80,20){\line(1,0){40}}\put(100,30){\msmash{a}}
    \put(120,20){\curve(0,0,20,3,40,0)\curve(0,0,20,-3,40,0)}\put(140,30){\msmash{b,c}}
    \put(160,20){\line(1,0){40}}\put(180,30){\msmash{a}}
    \put(200,20){\curve(0,0,20,3,40,0)\curve(0,0,20,-3,40,0)}\put(220,30){\msmash{b,d}}
    \put(240,20){\line(1,0){40}}\put(260,30){\msmash{a}}
    \put(0,10){\msmash{110}}
    \put(0,20){\spline(0,0)(15,25)(-15,25)(0,0)}\put(-25,35){\msmash{b,c,d}}
    \put(40,15){\msmash{010}}
    \put(40,20){\spline(0,0)(15,25)(-15,25)(0,0)}\put(40,35){\msmash{c}}
    \put(80,15){\msmash{000}}
    \put(80,20){\spline(0,0)(15,25)(-15,25)(0,0)}\put(80,35){\msmash{c}}
    \put(120,15){\msmash{100}}
    \put(120,20){\spline(0,0)(15,25)(-15,25)(0,0)}\put(120,35){\msmash{d}}
    \put(160,15){\msmash{101}}
    \put(160,20){\spline(0,0)(15,25)(-15,25)(0,0)}\put(160,35){\msmash{d}}
    \put(200,15){\msmash{001}}
    \put(200,20){\spline(0,0)(15,25)(-15,25)(0,0)}\put(200,35){\msmash{c}}
    \put(240,15){\msmash{011}}
    \put(240,20){\spline(0,0)(15,25)(-15,25)(0,0)}\put(240,35){\msmash{c}}
    \put(280,15){\msmash{111}}
    \put(280,20){\spline(0,0)(15,25)(-15,25)(0,0)}\put(302,35){\msmash{b,c,d}}
  \end{picture}
\end{center}  

The graphs $\gf_2$ and $\gf_3$ for $\tilde G$ are similar:
\begin{center}
  \begin{picture}(120,40)
    \put(0,20){\line(1,0){40}}\put(20,30){\msmash{a}}
    \put(40,20){\line(1,0){40}}\put(60,30){\msmash{b}}
    \put(80,20){\line(1,0){40}}\put(100,30){\msmash{a}}
    \put(0,15){\msmash{10}}
    \put(0,20){\spline(0,0)(18,27)(-18,27)(0,0)}\put(-25,35){\msmash{b,c,d}}
    \put(40,15){\msmash{00}}
    \put(40,20){\spline(0,0)(15,25)(-15,25)(0,0)}\put(40,38){\msmash{c,d}}
    \put(80,15){\msmash{01}}
    \put(80,20){\spline(0,0)(15,25)(-15,25)(0,0)}\put(80,38){\msmash{c,d}}
    \put(120,15){\msmash{11}}
    \put(120,20){\spline(0,0)(15,25)(-15,25)(0,0)}\put(142,35){\msmash{b,c,d}}
  \end{picture}\\
  \begin{picture}(280,40)
    \put(0,20){\line(1,0){40}}\put(20,30){\msmash{a}}
    \put(40,20){\line(1,0){40}}\put(60,30){\msmash{d}}
    \put(80,20){\line(1,0){40}}\put(100,30){\msmash{a}}
    \put(120,20){\line(1,0){40}}\put(140,30){\msmash{b}}
    \put(160,20){\line(1,0){40}}\put(180,30){\msmash{a}}
    \put(200,20){\line(1,0){40}}\put(220,30){\msmash{d}}
    \put(240,20){\line(1,0){40}}\put(260,30){\msmash{a}}
    \put(0,10){\msmash{110}}
    \put(0,20){\spline(0,0)(15,25)(-15,25)(0,0)}\put(-25,35){\msmash{b,c,d}}
    \put(40,15){\msmash{010}}
    \put(40,20){\spline(0,0)(15,25)(-15,25)(0,0)}\put(40,38){\msmash{b,c}}
    \put(80,15){\msmash{000}}
    \put(80,20){\spline(0,0)(15,25)(-15,25)(0,0)}\put(80,38){\msmash{b,c}}
    \put(120,15){\msmash{100}}
    \put(120,20){\spline(0,0)(15,25)(-15,25)(0,0)}\put(120,38){\msmash{c,d}}
    \put(160,15){\msmash{101}}
    \put(160,20){\spline(0,0)(15,25)(-15,25)(0,0)}\put(160,38){\msmash{c,d}}
    \put(200,15){\msmash{001}}
    \put(200,20){\spline(0,0)(15,25)(-15,25)(0,0)}\put(200,38){\msmash{b,c}}
    \put(240,15){\msmash{011}}
    \put(240,20){\spline(0,0)(15,25)(-15,25)(0,0)}\put(240,38){\msmash{b,c}}
    \put(280,15){\msmash{111}}
    \put(280,20){\spline(0,0)(15,25)(-15,25)(0,0)}\put(302,35){\msmash{b,c,d}}
  \end{picture}
\end{center}  

\subsection{$\sch(\Gamma,P,S)$}
Recall that for $\Gamma$ we take $\Sigma=\{0,1,2\}$. Let us write
$\gf_n=\sch(\Gamma,P_n,S)$. First, $\gf_0$ has one vertex, labeled by
the empty sequence $\emptyset$, and four loops, labeled
$a,a^{-1},t,t^{-1}$.  Next, $\gf_1$ has three vertices, labeled
$0,1,2$, cyclically connected by a triangle labeled $a,a^{-1}$, and
with two loops at each vertex, labeled $t,t^{-1}$. In the pictures
only geometrical edges, in pairs $\{a,a^{-1}\}$ and $\{t,t^{-1}\}$,
are represented.

Now given $\gf_n$, for some $n\ge1$, perform on it the following
transformation: replace the vertex-labels $\sigma$ by $2\sigma$;
replace all triangles labeled by $a,a^{-1}$ connecting
$\rho,\sigma,\tau$ by: three triangles labeled by $a,a^{-1}$
connecting respectively  $0\rho,1\rho,2\rho$ and
$0\sigma,1\sigma,2\sigma$ and $0\tau,1\tau,2\tau$; a triangle labeled
by $t,t^{-1}$ connecting $0\rho,0\sigma,0\tau$; and loops labeled by
$t,t^{-1}$ at $1\rho,1\sigma,1\tau$. We claim the resulting graph is
$\gf_{n+1}$.

As above, it suffices to check that the letters on the edge-labels act
as described on the vertex-labels. If $a(\rho)=\sigma$ and
$t(\rho)=\tau$, then $t(0\rho)=0\sigma$, $t(1\rho)=1\rho$ and
$t(2\sigma)=2\tau$. The verification for $a$ edges is even simpler.

\begin{figure}
  \epsfig{file=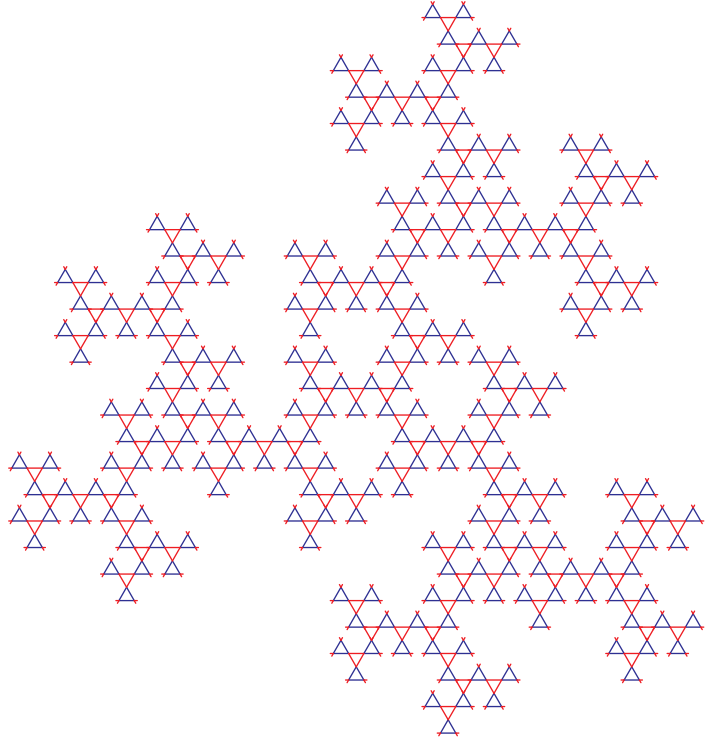}
  \caption{The Schreier Graph of $\Gamma_6$. The red and blue edges
    represent the generators \red{$s$} and \blue{$a$}.}
  \label{fig:schreier}
\end{figure}

\subsection{Substitutional graphs}
The three Schreier graphs presented in the previous subsection are
special cases of \emdef{substitutional graphs}, which we define below.

Substitutional graphs were introduced in the late 70's to describe
growth of multicellular organisms. They bear a strong similarity to
$L$-systems~\cite{rozenberg-s:l}, as was noted by Mikhael
Gromov~\cite{gromov:iggo}. Another notion of graph substitution has
been studied by~\cite{previte:subs}, where he had the same convergence
preoccupations as us.

Let us make a conventions in this section: all graphs
$\gf=(V(\gf),E(\gf))$ shall be \emdef{labeled}, i.e.\ endowed with
a map $E(\gf)\to C$ for a fixed set $C$ of colours, and
\emdef{pointed}, i.e.\ shall have a distinguished vertex $*\in
V(\gf)$. A \emdef{graph embedding} $\gf'\hookrightarrow\gf$ is just an
injective map $E(\gf')\to E(\gf)$ preserving the adjacency operations.
\begin{defn}
  A \emdef{substitutional rule} is a tuple $(U,R_1,\dots,R_n)$, where
  $U$ is a finite $d$-regular edge-labeled graph, called the
  \emdef{axiom}, and each $R_i$ is a rule of the form $X_i\to Y_i$,
  where $X_i$ and $Y_i$ are finite edge-labeled graphs. The graphs
  $X_i$ are required to have no common edge. Furthermore, there is a
  inclusion, written $\iota_i$, of the vertices of $X_i$ in the
  vertices of $Y_i$; the degree of $\iota_i(x)$ is the same as the
  degree of $x$ for all $x\in V(X_i)$, and all vertices of $Y_i$ not
  in the image of $\iota_i$ have degree $d$.
\end{defn}

Given a substitutional rule, one sets $\gf_0=U$ and constructs
iteratively $\gf_{n+1}$ from $\gf_n$ by listing all embeddings of all
$X_i$ in $\gf_n$ (they are disjoint), and replacing them by the
corresponding $Y_i$. If the base point $*$ of $\gf_n$ is in a graph $X_i$,
the base point of $\gf_{n+1}$ will be $\iota_i(*)$.

Note that this expansion operation preserves the degree, so $\gf_n$ is
a $d$-regular finite graph for all $n$. We are interested in fixed
points of this iterative process.

For any $R\in\N$, consider the balls $B_n(*,R)$ of radius $R$ at the
base point $*$ in $\gf_n$. Since there is only a finite number of
rules, there is an infinite sequence $n_0<n_1<\dots$ such that the
balls $B_{n_i}(*,R)\subseteq\gf_{n_i}$ are all isomorphic. We consider
$\gf$, a limit graph in the sense of~\cite{grigorchuk-z:infinite} (the
limit exists), and call it a \emph{substitutional graph}.

Note that in case some rule $X_i\to Y_i$ does not satisfy $X_i\subset
Y_i$, there may be more than one limit graph --- this is the case in
the following example (where there are in fact three limit graphs,
according to whether the next-to-rightmost edge is $b,c$ or $b,d$ or
$c,d$).
\begin{thm}
  The following four substitutional rules describe the Schreier graph of $G$:
  \begin{center}
    \begin{picture}(360,50)
      \put(0,40){\line(1,0){60}}\put(30,50){\msmash a}
      \put(0,45){\footnotesize\msmash\sigma}
      \put(60,45){\footnotesize\msmash\tau}
      \put(30,35){\vector(0,-1){20}}
      \put(80,40){\line(1,0){60}}\put(110,50){\msmash b}
      \put(80,45){\footnotesize\msmash\sigma}
      \put(140,45){\footnotesize\msmash\tau}
      \put(110,35){\vector(0,-1){20}}
      \put(160,40){\line(1,0){60}}\put(190,50){\msmash c}
      \put(160,45){\footnotesize\msmash\sigma}
      \put(220,45){\footnotesize\msmash\tau}
      \put(190,35){\vector(0,-1){20}}
      \put(240,40){\line(1,0){60}}\put(270,50){\msmash d}
      \put(240,45){\footnotesize\msmash\sigma}
      \put(300,45){\footnotesize\msmash\tau}
      \put(270,35){\vector(0,-1){20}}
      \put(0,10){\line(1,0){20}}\put(10,6){\msmash a}
      \put(20,10){\curve(0,0,10,2,20,0)\curve(0,0,10,-2,20,0)}
      \put(30,4){\msmash{b,c}}
      \put(40,10){\line(1,0){20}}\put(50,6){\msmash a}
      \put(20,10){\spline(0,0)(12,20)(-12,20)(0,0)}\put(10,20){\msmash d}
      \put(40,10){\spline(0,0)(12,20)(-12,20)(0,0)}\put(50,20){\msmash d}
      \put(0,7){\footnotesize\msmash{1\sigma}}
      \put(20,7){\footnotesize\msmash{0\sigma}}
      \put(40,7){\footnotesize\msmash{0\tau}}
      \put(60,7){\footnotesize\msmash{1\tau}}
      \put(80,10){\line(1,0){60}}\put(110,6){\msmash d}
      \put(80,7){\footnotesize\msmash{1\sigma}}
      \put(140,7){\footnotesize\msmash{1\tau}}
      \put(160,10){\line(1,0){60}}\put(190,6){\msmash b}
      \put(160,7){\footnotesize\msmash{1\sigma}}
      \put(220,7){\footnotesize\msmash{1\tau}}
      \put(240,10){\line(1,0){60}}\put(270,6){\msmash c}
      \put(240,7){\footnotesize\msmash{1\sigma}}
      \put(300,7){\footnotesize\msmash{1\tau}}
      \put(340,40){\msmash{\text{axiom}}}
      \put(320,10){\spline(0,0)(12,20)(-12,20)(0,0)}\put(340,22){\msmash{b,c,d}}
      \put(320,10){\line(1,0){40}}\put(340,6){\msmash a}
      \put(360,10){\spline(0,0)(12,20)(-12,20)(0,0)}
      \put(320,7){\footnotesize\msmash{0}}
      \put(360,7){\footnotesize\msmash{1}}
    \end{picture}
  \end{center}
  where the vertex inclusions are given by $\sigma\mapsto1\sigma$ and
  $\tau\mapsto1\tau$. The base point is the vertex $111\dots$.
\end{thm}
\begin{proof}
  Consider the Schreier graph $\gf_n$ associated to the action of $G$
  on $\Sigma^n$, the $n$-th level of the tree $\tree_\Sigma$. The
  vertex set of $\gf_n$ is $\Sigma^n$, and its edges are described by
  the action of $G$.  Note first that the axiom is $\gf_0$.
  
  We construct $\gf_{n+1}$ from $\gf_n$. Split $\Sigma^{n+1}$ as
  $0\Sigma^n\sqcup 1\Sigma^n$. By virtue of the definition of $\phi$
  given in~(\ref{eq:phiG}), the $b,c,d$-edges within $1\Sigma^n$ are
  in bijection to the $c,d,b$-edges in $\gf_n$, while the $b,c$-edges
  within $0\Sigma^n$ are in bijection with the $a$-edges in $\gf_n$,
  and there are $d$-loops at all $\sigma\in0\Sigma^n$. Moreover there
  are ``parallel edges'' labeled $a$ between $0\sigma$ and $1\sigma$
  for all $\sigma\in\Sigma^n$.
  
  Now consider any $b,c,d$-edge in $\gf_n$, say between $\sigma$ and
  $\tau$. In $\gf_{n+1}$, it gives rise to a $d,b,c$-edge between
  $1\sigma$ and $1\tau$.
  
  Consider then an $a$-edge in $\gf_n$ between $\sigma$ and $\tau$. In
  $\gf_{n+1}$, it gives rise to the following subgraph: an $a$-edge
  from $1\sigma$ to $0\sigma$; two edges, labeled $b$ and $c$, from
  $0\sigma$ to $0\tau$; an $a$-edge from $0\tau$ to $1\tau$; and two
  loops, labeled $d$ at $0\sigma$ and $0\tau$. This is precisely the
  substitutional rule for $a$, completing the proof.
\end{proof}

We omit the similar proof of the following result:
\begin{thm}
  The following four substitutional rules describe the Schreier graph
  of $\overline G$:
  \begin{center}
    \begin{picture}(360,50)
      \put(0,40){\line(1,0){60}}\put(30,50){\msmash a}
      \put(0,45){\footnotesize\msmash\sigma}
      \put(60,45){\footnotesize\msmash\tau}
      \put(30,35){\vector(0,-1){20}}
      \put(80,40){\line(1,0){60}}\put(110,50){\msmash{\tilde b}}
      \put(80,45){\footnotesize\msmash\sigma}
      \put(140,45){\footnotesize\msmash\tau}
      \put(110,35){\vector(0,-1){20}}
      \put(160,40){\line(1,0){60}}\put(190,50){\msmash{\tilde c}}
      \put(160,45){\footnotesize\msmash\sigma}
      \put(220,45){\footnotesize\msmash\tau}
      \put(190,35){\vector(0,-1){20}}
      \put(240,40){\line(1,0){60}}\put(270,50){\msmash{\tilde d}}
      \put(240,45){\footnotesize\msmash\sigma}
      \put(300,45){\footnotesize\msmash\tau}
      \put(270,35){\vector(0,-1){20}}
      \put(0,10){\line(1,0){20}}\put(10,6){\msmash a}
      \put(20,10){\line(1,0){20}}\put(30,4){\msmash{\tilde b}}
      \put(40,10){\line(1,0){20}}\put(50,6){\msmash a}
      \put(20,10){\spline(0,0)(12,20)(-12,20)(0,0)}\put(2,25){\msmash{\tilde c,\tilde d}}
      \put(40,10){\spline(0,0)(12,20)(-12,20)(0,0)}\put(58,25){\msmash{\tilde c,\tilde d}}
      \put(0,7){\footnotesize\msmash{1\sigma}}
      \put(20,7){\footnotesize\msmash{0\sigma}}
      \put(40,7){\footnotesize\msmash{0\tau}}
      \put(60,7){\footnotesize\msmash{1\tau}}
      \put(80,10){\line(1,0){60}}\put(110,4){\msmash{\tilde d}}
      \put(80,7){\footnotesize\msmash{1\sigma}}
      \put(140,7){\footnotesize\msmash{1\tau}}
      \put(160,10){\line(1,0){60}}\put(190,4){\msmash{\tilde b}}
      \put(160,7){\footnotesize\msmash{1\sigma}}
      \put(220,7){\footnotesize\msmash{1\tau}}
      \put(240,10){\line(1,0){60}}\put(270,6){\msmash{\tilde c}}
      \put(240,7){\footnotesize\msmash{1\sigma}}
      \put(300,7){\footnotesize\msmash{1\tau}}
      \put(340,40){\msmash{\text{axiom}}}
      \put(320,10){\spline(0,0)(12,20)(-12,20)(0,0)}\put(340,22){\msmash{\tilde b,\tilde c,\tilde d}}
      \put(320,10){\line(1,0){40}}\put(340,6){\msmash a}
      \put(360,10){\spline(0,0)(12,20)(-12,20)(0,0)}
      \put(320,7){\footnotesize\msmash{0}}
      \put(360,7){\footnotesize\msmash{1}}
    \end{picture}
  \end{center}
  where the vertex inclusions are given by $\sigma\mapsto1\sigma$ and
  $\tau\mapsto1\tau$. The base point is the vertex $111\dots$.
\end{thm}

\begin{thm}
  The substitutional rules producing the Schreier graphs of $\Gamma$
  and $\overline\Gamma$ are given below. Remember that the Schreier
  graphs of $\overline\Gamma$ and $\doverline\Gamma$ are isomorphic:
  \begin{center}
    \begin{picture}(300,100)(0,-10)
      \put(0,40){\msmash{\Gamma:}}
      \put(0,0){\blue{\line(1,0){100}}}\put(50,-4){\blue{\msmash a}}
      \put(0,0){\blue{\line(15,26){50}}}\put(15,43){\blue{\msmash a}}
      \put(100,0){\blue{\line(-15,26){50}}}\put(85,43){\blue{\msmash a}}
      \put(0,-3){\footnotesize\msmash\rho}
      \put(100,-3){\footnotesize\msmash\sigma}
      \put(50,92){\footnotesize\msmash\tau}
      \put(100,40){\vector(1,0){40}}
      \put(130,0){\blue{\line(1,0){40}}}\put(150,7){\blue{\msmash a}}
      \put(130,0){\blue{\line(15,26){20}}}\put(143,17){\blue{\msmash a}}
      \put(170,0){\blue{\line(-15,26){20}}}\put(157,17){\blue{\msmash a}}
      \put(130,-3){\footnotesize\msmash{2\rho}}
      \put(144,37){\footnotesize\msmash{1\rho}}
      \put(168,-3){\footnotesize\msmash{0\rho}}
      \put(190,0){\blue{\line(1,0){40}}}\put(210,7){\blue{\msmash a}}
      \put(190,0){\blue{\line(15,26){20}}}\put(203,17){\blue{\msmash a}}
      \put(230,0){\blue{\line(-15,26){20}}}\put(217,17){\blue{\msmash a}}
      \put(195,-3){\footnotesize\msmash{1\sigma}}
      \put(217,37){\footnotesize\msmash{0\sigma}}
      \put(230,-3){\footnotesize\msmash{2\sigma}}
      \put(160,52){\blue{\line(1,0){40}}}\put(180,59){\blue{\msmash a}}
      \put(160,52){\blue{\line(15,26){20}}}\put(173,69){\blue{\msmash a}}
      \put(200,52){\blue{\line(-15,26){20}}}\put(187,69){\blue{\msmash a}}
      \put(158,57){\footnotesize\msmash{0\tau}}
      \put(180,92){\footnotesize\msmash{2\tau}}
      \put(205,57){\footnotesize\msmash{1\tau}}
      \put(170,0){\red{\line(-10,52){10}}}\put(157,37){\red{\msmash s}}
      \put(170,0){\red{\line(40,35){40}}}\put(185,7){\red{\msmash s}}
      \put(160,52){\red{\line(50,-17){50}}}\put(195,48){\red{\msmash s}}
      \put(150,35){\red{\spline(0,0)(10,17)(-10,17)(0,0)}}
      \put(190,0){\red{\spline(0,0)(-20,0)(-10,-17)(0,0)}}
      \put(200,52){\red{\spline(0,0)(20,0)(10,-17)(0,0)}}
      \put(300,70){\msmash{\text{axiom}}}
      \put(280,0){\blue{\line(1,0){40}}}\put(300,7){\blue{\msmash a}}
      \put(280,0){\blue{\line(15,26){20}}}\put(293,17){\blue{\msmash a}}
      \put(320,0){\blue{\line(-15,26){20}}}\put(307,17){\blue{\msmash a}}
      \put(283,-3){\footnotesize\msmash{2}}\put(272,7){\red{\msmash s}}
      \put(294,37){\footnotesize\msmash{1}}\put(328,7){\red{\msmash s}}
      \put(318,-3){\footnotesize\msmash{0}}\put(310,45){\red{\msmash s}}
      \put(300,35){\red{\spline(0,0)(10,17)(-10,17)(0,0)}}
      \put(280,0){\red{\spline(0,0)(-20,0)(-10,-17)(0,0)}}
      \put(320,0){\red{\spline(0,0)(20,0)(10,-17)(0,0)}}
    \end{picture}\\[2mm]
    \begin{picture}(300,100)(0,-10)
      \put(0,40){\msmash{\overline\Gamma,\doverline\Gamma:}}
      \put(0,0){\blue{\line(1,0){100}}}\put(50,-4){\blue{\msmash a}}
      \put(0,0){\blue{\line(15,26){50}}}\put(15,43){\blue{\msmash a}}
      \put(100,0){\blue{\line(-15,26){50}}}\put(85,43){\blue{\msmash a}}
      \put(0,-3){\footnotesize\msmash\rho}
      \put(100,-3){\footnotesize\msmash\sigma}
      \put(55,87){\footnotesize\msmash\tau}
      \put(100,40){\vector(1,0){40}}
      \put(130,0){\blue{\line(1,0){40}}}\put(150,7){\blue{\msmash a}}
      \put(130,0){\blue{\line(15,26){20}}}\put(143,17){\blue{\msmash a}}
      \put(170,0){\blue{\line(-15,26){20}}}\put(157,17){\blue{\msmash a}}
      \put(130,-3){\footnotesize\msmash{2\rho}}
      \put(144,37){\footnotesize\msmash{1\rho}}
      \put(168,-3){\footnotesize\msmash{0\rho}}
      \put(190,0){\blue{\line(1,0){40}}}\put(210,7){\blue{\msmash a}}
      \put(190,0){\blue{\line(15,26){20}}}\put(203,17){\blue{\msmash a}}
      \put(230,0){\blue{\line(-15,26){20}}}\put(217,17){\blue{\msmash a}}
      \put(195,-3){\footnotesize\msmash{1\sigma}}
      \put(217,37){\footnotesize\msmash{0\sigma}}
      \put(230,-3){\footnotesize\msmash{2\sigma}}
      \put(160,52){\blue{\line(1,0){40}}}\put(180,59){\blue{\msmash a}}
      \put(160,52){\blue{\line(15,26){20}}}\put(173,69){\blue{\msmash a}}
      \put(200,52){\blue{\line(-15,26){20}}}\put(187,69){\blue{\msmash a}}
      \put(158,57){\footnotesize\msmash{0\tau}}
      \put(186,87){\footnotesize\msmash{2\tau}}
      \put(205,57){\footnotesize\msmash{1\tau}}
      \put(170,0){\red{\line(-10,52){10}}}\put(159,35){\red{\msmash t}}
      \put(190,0){\red{\line(10,52){10}}}
      \put(170,0){\red{\line(40,35){40}}}\put(188,10){\red{\msmash t}}
      \put(190,0){\red{\line(-40,35){40}}}
      \put(160,52){\red{\line(50,-17){50}}}\put(194,46){\red{\msmash t}}
      \put(200,52){\red{\line(-50,-17){50}}}
      \put(300,70){\msmash{\text{axiom}}}
      \put(280,0){\blue{\line(1,0){40}}}\put(300,7){\blue{\msmash a}}
      \put(280,0){\blue{\line(15,26){20}}}\put(293,17){\blue{\msmash a}}
      \put(320,0){\blue{\line(-15,26){20}}}\put(307,17){\blue{\msmash a}}
      \put(283,-3){\footnotesize\msmash{2}}\put(272,7){\red{\msmash t}}
      \put(294,37){\footnotesize\msmash{1}}\put(328,7){\red{\msmash t}}
      \put(318,-3){\footnotesize\msmash{0}}\put(310,45){\red{\msmash t}}
      \put(300,35){\red{\spline(0,0)(10,17)(-10,17)(0,0)}}
      \put(280,0){\red{\spline(0,0)(-20,0)(-10,-17)(0,0)}}
      \put(320,0){\red{\spline(0,0)(20,0)(10,-17)(0,0)}}
    \end{picture}
  \end{center}
  where the vertex inclusions are given by $\rho\mapsto2\rho$,
  $\sigma\mapsto2\sigma$ and $\tau\mapsto2\tau$. The base point is the
  vertex $222\dots$.
\end{thm}
\begin{proof}
  We prove the claim for $\Gamma$ only, as the same reasoning applies
  for $\overline\Gamma$. Consider the Schreier graph $\gf_n$
  associated to the action of $G$ on $\Sigma^n$, the $n$-th level of
  the tree $\tree_\Sigma$. The vertex set of $\gf_n$ is $\Sigma^n$,
  and its edges are described by the action of $G$.  Note first that
  the axiom is $\gf_0$.
  
  We construct $\gf_{n+1}$ from $\gf_n$. Split $\Sigma^{n+1}$ as
  $0\Sigma^n\sqcup 1\Sigma^n\sqcup2\Sigma^n$. By virtue of the
  definition of $\phi$ given in~(\ref{eq:phiGamma}), the $s$-edges
  within $2\Sigma^n$ are in bijection to the $s$-edges in $\gf_n$,
  while the $s$-edges within $0\Sigma^n$ are in bijection with the
  $a$-edges in $\gf_n$, and there are $s$-loops at every
  $\sigma\in1\Sigma^n$. Moreover there are ``parallel triangles''
  labeled $a$ between $0\sigma$, $1\sigma$ and $2\sigma$ for all
  $\sigma\in\Sigma^n$.

  Now consider any $s$-edge in $\gf_n$, say between $\sigma$ and
  $\tau$. In $\gf_{n+1}$, it remains an $s$-edge, but now between
  $2\sigma$ and $2\tau$.
  
  Consider then an $a$-edge in $\gf_n$ between $\sigma$ and $\tau$. In
  $\gf_{n+1}$, it gives rise to the following subgraph: an
  $a$-triangle between $0\sigma$, $1\sigma$ and $2\sigma$; an $s$-edge
  between $0\sigma$ and $0\tau$; $s$-loops at $1\sigma$ and $1\tau$;
  and an $a$-triangle between $0\tau$, $1\tau$ and $2\tau$. Actually
  the $a$-edges form triangles so these subgraphs overlap at the
  $a$-triangles and $s$-loops. This justifies the substitutional rule
  for $a$-triangles, completing the proof.
\end{proof}

By Proposition~\ref{prop:polygrowth}, the limit graphs have
asymptotically polynomial growth of degree $\log_2(3)$.

Note that there are maps $\pi_n:V(\gf_{n+1})\to V(\gf_n)$ that locally
(i.e.\ in each copy of some right-hand rule $Y_i$) are the inverse of
the embedding $\iota_i$. In case these $\pi_n$ are graph morphisms,
and one can consider the projective system $\{\gf_n,\pi_n\}$, and its
inverse limit $\widehat\gf=\varprojlim\gf_n$, which is a profinite
graph~\cite{ribes-z:profinite}. We devote our attention to the
discrete graph $\gf=\varinjlim\gf_n$.

The growth series of $\gf$ can often be described as an infinite
product. We give such an expression for the graph in
Figure~\ref{fig:schreier}, making use of the fact that $\gf$ ``looks
like a tree'' (even though it is amenable).

Consider the finite graphs $\gf_n$; recall that $\gf_n$ has $3^n$
vertices. Let $D_n$ be the diameter of $\gf_n$ (maximal distance
between two vertices), and let $\gamma_n=\sum_{i\in\N}\gamma_n(i)X^i$
be the growth series of $\gf_n$ (here $\gamma_n(i)$ denotes the number
of vertices in $\gf_n$ at distance $i$ from the base point $*$).

The construction rule for $\gf$ implies that $\gf_{n+1}$ can be
constructed as follows: take three copies of $\gf_n$, and in each of
them mark a vertex $V$ at distance $D_n$ from $*$. At each $V$ delete
the loop labeled $s$, and connect the three copies by a triangle
labeled $s$ at the three $V$'s.
It then follows that $D_{n+1}=2D_n+1$, and
$\gamma_{n+1}=(1+2X^{D_n+1})\gamma_n$. Using the initial values
$\gamma_0=1$ and $D_0=0$, we obtain by induction
\[D_n=2^n-1,\qquad\gamma_n=\prod_{i=0}^{n-1}(1+2X^{2^i}).\]
We also
have shown that the ball of radius $2^n$ around $*$ contains $3^n$
points, so the growth of $\gf$ is at least $n^{log_2(3)}$. But
Proposition~\ref{prop:polygrowth} shows that it is also an upper
bound, and we conclude:
\begin{prop}
  $\Gamma$ is an amenable $4$-regular graph whose growth function is
  transcendental, and admits the product decomposition
  \[\gamma(X)=\prod_{i\in\N}(1+2X^{2^i}).\]
  It is planar, and has polynomial growth of degree $\log_2(3)$.
\end{prop}

Any graph is a metric space when one identifies each edge with a
disjoint copy of an interval $[0,L]$ for some $L>0$. We turn $\gf_n$
in a diameter-$1$ metric space by giving to each edge in $\gf_n$ the
length $L=\operatorname{diam}(\gf_n)^{-1}$. The family $\{\gf_n\}$
then converges, in the following sense:

Let $A,B$ be closed subsets of the metric space $(X,d)$. For any
$\epsilon$, let $A_\epsilon=\{x\in X|\,d(x,A)\le\epsilon\}$, and
define the \emdef{Hausdorff distance}
\[d_X(A,B) = \inf\{\epsilon|\,A\subseteq B_\epsilon,B\subseteq
A_\epsilon\}.\]
This defines a metric on closed subsets of $X$. For
general metric spaces $(A,d)$ and $(B,d)$, define their
\emdef{Gromov-Hausdorff distance}
\[d^{GH}(A,B) = \inf_{X,i,j}d_X(i(A),j(B)),\]
where $i$ and $j$ are isometric embeddings of $A$ and $B$ in a metric
space $X$.

We may now rephrase the considerations above as follows: the sequence
$\{\gf_n\}$ is convergent in the Gromov-Hausdorff metric. The limit
set $\gf_\infty$ is a compact metric space.

The limit spaces are then: for $G$ and $\tilde G$, the limit
$\gf_\infty$ is the interval $[0,1]$ (in accordance with its linear
growth, see Proposition~\ref{prop:polygrowth}). The limit spaces for
$\Gamma$, $\overline\Gamma$ and $\doverline\Gamma$ are fractal sets of
dimension $\log_2(3)$.

When the present work was completed Stanislav Smirnov informed us
about the article~\cite{malozemov:fractal} which has some connection
to our paper. A substitutional tree $\tree$ is constructed
in~\cite{malozemov:fractal}, and the spectrum of the Markov operator
on $\tree$ is of the form $K\sqcup P$, where $K$ is the Julia set of
some quadratic map and $P$ is countable. The approximation arguments
used in~\cite{malozemov:fractal} are different from ours and does not
use any amenability assumption (indeed, the graphs constructed there
have eigenvalues with compactly supported eigenvectors).

An essential difference is that the graph in~\cite{malozemov:fractal}
is not regular and therefore is not a Schreier graph of a group.

\section{Concluding Remarks and Problems}
Let us draw here some conclusions and formulate a few questions. We
hope that our methods will prove useful in the investigation of
spectral properties of the Laplace operator $\Delta$ for other groups acting
on rooted trees. By computing the spectrum of the Hecke type operator of a
quasi-regular representation $\rho_{G/P}$ we obtained a subset of the
spectrum of $\Delta$; namely,
$\spec(\rho_{G/P})\subseteq\spec(\rho_G)$ as soon as $P$ is amenable.
Moreover, in one case (that of $\tilde G$) we obtained $\spec(\Delta)$
using the bipartitivity of $\tilde G$'s Cayley graph. More research
must be done in the nonamenable case, in particular
\begin{question} Under which conditions (besides amenability) does one
  have $\spec(\rho_{G/P})\subseteq\spec(\rho_G)$? When does one
  have $\spec(\rho_{G/P})=\spec(\rho_G)$?
\end{question}

We also hope that the methods in this article can be useful to find
residually finite examples of nonamenable groups without free
subgroups (all known examples are non-residually-finite).  The
non-amenability of the group considered could be proven by showing
that $1$ is an isolated point in the spectrum of the associated
dynamical system.
\begin{question} Is there a finitely generated subgroup $G$ of
  $\aut(\tree)$ with no non-abelian free subgroup and such that
  $1\notin\spec(\rho_{G/P})$?
\end{question}

Examples of groups $G$ with $1\notin\spec(\rho_{G/P})$ are also
interesting because they provide sequences of expanding
graphs~\cite{lubotzky:expanders}, namely the Schreier graphs
$\sch(G,P_n,S)$. In some cases these graphs can be Ramanujan graphs.
We formulate the following question with Andrzej \.Zuk:
\begin{question}
  Are there groups $G<\aut(\tree)$, generated by a finite set $S$ of
  finite automata, such that the sequence of graphs $\sch(G,P_n,S)$
  is \emph{(i)} a sequence of expanders; \emph{(ii)} a sequence of
  Ramanujan graphs?
\end{question}

The problem of factoring the resolvent of operators of Hecke type is
related to the decomposition of $\rho_{G/P_n}$ in irreducible
representations. It can be shown that any finite-dimensional
irreducible representation of $G$ is a subrepresentation in a tensor
product of sufficiently many copies of $\rho_{G/P_n}$ (in other words,
the representation ring $R(G)$ is generated by the irreducible
subrepresentations of $\rho_{G/P_n}$). We therefore hope that our
methods will extend the knowledge on the irreducible representations
of $G$.

We constructed a virtually torsion-free group, $\Gamma$, with totally
disconnected spectrum $\spec(\rho_{\Gamma/P})$. To answer in the
negative to the Baum-Connes conjecture, as well as to the
Kadison-Kaplansky conjecture~\cite{valette:idempotents}, it would
suffice to construct a torsion-free group with a gap in the spectrum
of its Laplace operator. Therefore one should delete the ``virtually''
and replace $\rho_{\Gamma/P}$ by $\rho_\Gamma$ in order to
improve our results.
\begin{question}
  Is there a torsion-free group $G<\aut(\tree)$ with totally
  disconnected spectrum $\spec(\rho_G)$? Or with a gap in its
  spectrum?
\end{question}

%

\section{Acknowledgments}
Both authors wish to thank heartily Pierre de la Harpe for his
constant availability and interest. Viviane Baladi, Marc Burger,
Vaughan Jones, Volodymyr Nekrashevych, Gilles Robert, Alain Valette
and Andrzej \.Zuk contributed by generous discussions and input.

\bibliographystyle{amsalpha}
\bibliography{../mrabbrev,../math}
\end{document}